\documentclass[12pt,a4paper,dvipdfm]{article}
\usepackage{amssymb,amsfonts}

\usepackage{amsmath, latexsym}

\usepackage{url}
\usepackage{gastex}

\newenvironment{Proof}{\rm \trivlist \item[\hskip \labelsep{\bf
Proof.}]}{\cqfd\endtrivlist}

\def\cqfd{\skip10=\parfillskip\parfillskip=0pt
\enspace\hfill\symbolecqfd\par\parfillskip=\skip10\par\medskip}
\def\symbolecqfd{\rlap{$\sqcap$}$\sqcup$}
\def\preuve{\begin{Proof}}
\def\proof{\begin{Proof}}
\def\eop{\end{Proof}}

\let\demo\proof
\let\qed\eop

\newtheorem{thm}{Theorem}[section]
\newtheorem{prop}[thm]{Proposition}
\newtheorem{lem}[thm]{Lemma}
\newtheorem{cor}[thm]{Corollary}
\newtheorem{obs}[thm]{Fact}

\newtheorem{pro-example}[thm]{Example}
\newenvironment{example}{\begin{pro-example}\rm}{\cqfd\end{pro-example}}

\newtheorem{pro-remark}[thm]{Remark}
\newenvironment{remark}{\begin{pro-remark}\rm}{\cqfd\end{pro-remark}}

\newcommand{\p}{\mathcal{P}}

\newcommand{\calE}{\mathcal{E}}
\newcommand{\calL}{\mathcal{L}}
\newcommand{\calO}{\mathcal{O}}
\newcommand{\calP}{\mathcal{P}}
\newcommand{\calT}{\mathcal{T}}
\let\phi\varphi
\def\leff{\le_{\textsf{ff}}}
\def\lefg{\le_{\textsf{fg}}}
\def\lefi{\le_{\textsf{fi}}}
\def\lealg{\le_{\textsf{alg}}}
\def\leealg{\le_{\textsf{ealg}}}
\let\lalg\lealg
\let\lff\leff
\let\lfi\lefi

\def\AE{\textsf{AE}}

\def\rank{\textsf{rk}}
\def\Aut{\textsf{Aut}}

\def\V{\textbf{V}}

\def\cl{\textsf{cl}}
\def\id{\textsf{id}}

\def\inv{^{-1}}

\def\note#1{}

\begin{document}

\title{Algebraic extensions in free groups}

\markright{\protect Algebraic extensions in free groups}

\author{Alexei Miasnikov
\small{(\url{alexeim@math.mcgill.ca})}\protect\footnote{Dept of
Mathematics and Statistics -- McGill University -- Burnside Hall, room
915 --  805 Sherbrooke West -- Montr\'eal, Quebec, Canada, H3A 2K6} \\
\small{Dept of Mathematics and Statistics, McGill University and}\\
\small{Dept of Mathematics and Computer Science, City University of
New York}\\
\\ Enric Ventura
\small{(\url{enric.ventura@upc.edu})}\protect\footnote{EPSEM,
Universitat Polit\`ecnica de Catalunya -- Av.  Bases de Manresa 61 --
73 08242 Manresa, Barcelona -- Spain} \\
\small{Universitat Polit\`ecnica de Catalunya and Centre de
Recerca Matem\`atica}\\
\\ Pascal Weil
\small{(\url{pascal.weil@labri.fr})}\protect\footnote{LaBRI -- 351
cours de la Lib\'eration -- 33405 Talence Cedex -- France} \\
\small{LaBRI (CNRS and Universit\'e Bordeaux-1)}}

\date{}

\maketitle

\begin{abstract}
The aim of this paper is to unify the points of view of three recent
and independent papers (Ventura 1997, Margolis, Sapir and Weil 2001
and Kapovich and Miasnikov 2002), where similar modern versions of a
1951 theorem of Takahasi were given.  We develop a theory of algebraic
extensions for free groups, highlighting the analogies and differences
with respect to the corresponding classical field-theoretic notions,
and we discuss in detail the notion of algebraic closure.  We apply
that theory to the study and the computation of certain algebraic
properties of subgroups (\textit{e.g.} being malnormal, pure, inert or
compressed, being closed in certain profinite topologies) and the
corresponding closure operators.  We also analyze the closure of a
subgroup under the addition of solutions of certain sets of equations.
\end{abstract}

\section{Introduction}

A well-known result by Nielsen and Schreier states that all subgroups
of a free group $F$ are free.  A non-specialist in group theory could
be tempted to guess from this pleasant result that the lattice of
subgroups of $F$ is simple, and easy to understand.  This is however
far from being the case, and a closer look quickly reveals the
classical fact that inclusions do not respect rank.  In fact, the free
group of countably infinite rank appears many times as a subgroup of
the free group of rank 2.  There are also many examples of subgroups
$H, K$ of $F$ such that the rank of $H\cap K$ is greater than the
ranks of $H$ and $K$.  These are just a few indications that the
lattice of subgroups of $F$ is not easy.

Although the lattice of subgroups of free groups was already studied
by earlier authors, Serre and Stallings in their seminal 
1977 and 1983 papers~\cite{Serre,St},
introduced a powerful new technique, that has since turned
out to be extremely useful in this line of research.  It consists in
thinking of $F$ as the fundamental group of a bouquet of circles $R$,
and of subgroups of $F$ as covering spaces of $R$, i.e. some special
types of graphs.  With this idea in mind, one can understand and prove
many properties of the lattice of subgroups of $F$ using graph theory.
These techniques are also very useful to solve algorithmic problems
and to effectively compute invariants concerning subgroups of $F$.

The present paper offers a contribution in this direction, by
analyzing a tool (an invariant associated to a given subgroup $H\leq
F$) which is suggested by a 1951 theorem of Takahasi~\cite{T} (see
Section~\ref{Takahasi section}).  The algorithmic constructions
involved in the computation of this invariant actually appeared in
recent years, in three completely independent papers~\cite{V},
\cite{MSW} and~\cite{KM}, where the same notion was invented in
independent ways.  In chronological order, we refer:
\begin{itemize}
\item to the \emph{fringe of a subgroup}, constructed in 1997 by 
Ventura (see~\cite{V}), and applied to the study of maximal rank fixed
subgroups of automorphisms of free groups;

\item to the \emph{overgroups} of a subgroup, constructed in 2001 by
Margolis, Sapir and Weil (see~\cite{MSW}), and applied to
improve an algorithm of Ribes-Zaleskii for computing the pro-$p$
topological closure of a finitely generated subgroup of a free group,
among other applications; and

\item to the \textit{algebraic extensions} constructed in 2002 by
Kapovich and Miasnikov (see~\cite{KM}), in the context of a paper
where the authors surveyed, clarified and extended the list of
Stallings graphical techniques.
\end{itemize}

Turner also used the same notion, restricted to the case of
cyclic subgroups, in his paper~\cite{Tu} (again, independently) when
trying to find examples of test elements for the free group.

The terminology and the notation used in the above mentioned papers
are different, but the basic concept -- that of algebraic extension
for free groups -- is the same.  Although aimed at different
applications, the underlying basic result in these three papers is a
modern version of an old theorem by Takahasi~\cite{T}.  It states
that, for every finitely generated subgroup $H$ of a free group $F$,
there exist finitely many subgroups $H_0, \ldots , H_n$ canonically
associated to $H$, such that every subgroup of $F$ containing $H$ is a
free multiple of $H_i$ for some $i=0,\ldots,n$.  The original proof
was combinatorial, while the proof provided in~\cite{V}, \cite{MSW}
and~\cite{KM} (which is the same up to technical details) is
graphical, algorithmic, simpler and more natural.

The aim of this paper is to unify the points of view
in~\cite{V}, \cite{MSW} and~\cite{KM}, and to systematize the study of
the concept of algebraic extensions in free groups.  We show how
algebraic extensions intervene in the computation of certain abstract
closure properties for subgroups, sometimes making these properties
decidable.  This was the idea behind the application of algebraic
extensions to the study of profinite topological closures
in~\cite{MSW}, but it can be applied in other contexts.  In
particular, we extend the discussion of the notions of pure closure,
malnormal closure, inert closure, etc (a discussion that was initiated
in~\cite{KM}).

A particularly interesting application concerns the property of being
clo\-sed under the addition of the solutions of certain sets of
equations.  In this case, new results are obtained, and in particular
one can show that the rank of the closure of a subgroup $H$ is at most
equal to $\rank(H)$.

The paper is organized as follows.

In section~\ref{sec preliminaries}, we remind the readers of the
fundamentals of the representation of finitely generated subgroups of
a free group $F$ by finite labeled graphs.  This method, which was
initiated by Serre and Stallings at the end of the 1970s, quickly
became one of the major tools of the combinatorial theory of free
groups.  This leads us to the short, algorithmic proof of Takahasi's
theorem discussed above (see Section~\ref{Takahasi section}).

Section~\ref{sec algebraic extensions} introduces algebraic
extensions, essentially as follows: the algebraic extensions of a
finitely generated subgroup $H$ are the minimum family that can be
associated to $H$ by Takahasi's theorem.  We also discuss the
analogies that arise between this notion of algebraic extensions and
classical field-theoretic notions, and we discuss in detail the
corresponding notion of algebraic closure.

Section~\ref{sec abstract} is devoted to the applications of algebraic
extensions.  We show that whenever an abstract property of subgroups
of free groups is closed under free products and finite intersections,
then every finitely generated subgroup $H$ admits a unique closure
with respect to this property, which is finitely generated and is one
of the algebraic extensions of $H$.  Examples of such properties
include malnormality, purity or inertness, as well as the property of
being closed for certain profinite topologies.  In a number of
interesting situations, this leads to simple decidability results.
Equations over a subgroup, or rather the property of being closed
under the addition of solutions of certain sets of equations, provide
another interesting example of such an abstract property of subgroups,
which we discuss in Section~\ref{sec equations}.

Finally, in section~\ref{sec open}, we collect the open questions and
conjectures suggested by previous sections.

\section{Preliminaries}\label{sec preliminaries}

Throughout this paper, $A$ is a finite non-empty set and $F(A)$ (or
simply $F$ if no confusion may arise) is the free group on $A$.

In the algorithmic or computational statements on subgroups of free
groups, we tacitly assume that the free group $F$ is given together
with a basis $A$, that the elements of $F$ are expressed as words over
$A$, and that finitely generated subgroups of $F$ are given to us by
finite sets of generators, and hence by finite sets of words.

\subsection{Representation of subgroups of free
groups}\label{geometrical apparatus}

In his 1983 paper \cite{St}, Stallings showed how many of the
algorithmic constructions introduced in the first half of the 20th
century to handle finitely generated subgroups of free groups, can be
clarified and simplified by adopting a graph-theoretic language.  This
method has been used since then in a vast array of articles, including
work by the co-authors of this paper.

The fundamental notion is the existence of a natural, algorithmically
simple one-to-one correspondence between subgroups of the free group
$F$ with basis $A$, and certain $A$-labeled graphs -- mapping finitely
generated subgroups to finite graphs and vice versa.  This is nothing
else than a particular case of the more general covering theory for
topological spaces, particularized to graphs and free groups.  We
briefly describe this correspondence in the rest of this subsection.
More detailed expositions can be found in the literature: see
Stallings~\cite{St} or~\cite{V,KM} for a graph-oriented version, and
see one of~\cite{BMMW,Weil00,MSW,SW} for a more combinatorial-oriented
version, written in the language of automata theory.

By an \emph{$A$-labeled graph} $\Gamma$ we understand a directed graph
(allowing loops and multiple edges) with a designated vertex written
1, and in which each edge is labeled by a letter of $A$.  We say that
$\Gamma$ is \emph{reduced} if it is connected (more precisely, the
underlying undirected graph is connected), if distinct edges with the
same origin (resp.  with the same end vertex) have distinct labels,
and if every vertex $v\ne 1$ is adjacent to at least two different
edges.

In an $A$-labeled graph, we consider paths, where we are allowed to
travel backwards along edges.  The \textit{label} of such a path $p$
is the word obtained by reading consecutively the labels of the edges
crossed by $p$, reading $a^{-1}$ whenever an edge labeled $a\in A$ is
crossed backwards.  The path $p$ is called \textit{reduced} if it does
not cross twice consecutively the same edge, once in one direction and
then in the other.  Note that if $\Gamma$ is reduced then every
reduced path labels a reduced word in $F(A)$.

The subgroup of $F(A)$ associated with a reduced $A$-labeled graph
$\Gamma$ is the set of (reduced) words, which label reduced paths in
$\Gamma$ from the designated vertex 1 back to itself.  One can show
that every subgroup of $F(A)$ arises in this fashion, in a unique way.
That is, for each subgroup $H$ of $F(A)$, there exists a unique
reduced $A$-labeled graph, written $\Gamma_{A}(H)$, whose set of
labels of reduced closed paths at 1 is exactly $H$.

Moreover, if the subgroup $H$ is given together with a finite set of
generators $\{h_1, \ldots , h_r \}$ (where the $h_{i}$ are non-empty
reduced words over the alphabet $A\sqcup A\inv$), then one can
effectively construct $\Gamma_{A}(H)$, proceeding as follows.  First,
one constructs $r$ subdivided circles around a common distinguished
vertex 1, each labeled by one of the $h_{i}$ (and following the above
convention: an inverse letter, say $a\inv$ with $a\in A$, in a word
$h_{i}$ gives rise to an $a$-labeled edge in the reverse direction on
the corresponding circle).  If $h_{i}$ has length $n_{i}$, then the
corresponding circle has $n_{i}$ edges and $n_{i}-1$ vertices, in
addition to the vertex 1.  Then, we iteratively identify identically
labeled pairs of edges starting (resp.  ending) at the same vertex.
One shows that this process terminates, that it does not matter in
which order identifications take place, and that the resulting
$A$-labeled graph is reduced and equal to $\Gamma_{A}(H)$.  In
particular, it does not depend on the choice of a set of generators of
$H$.  Also, this shows that $\Gamma_{A}(H)$ is finite if and only if
$H$ is finitely generated (see one
of~\cite{St,V,BMMW,Weil00,MSW,KM,SW} for more details).

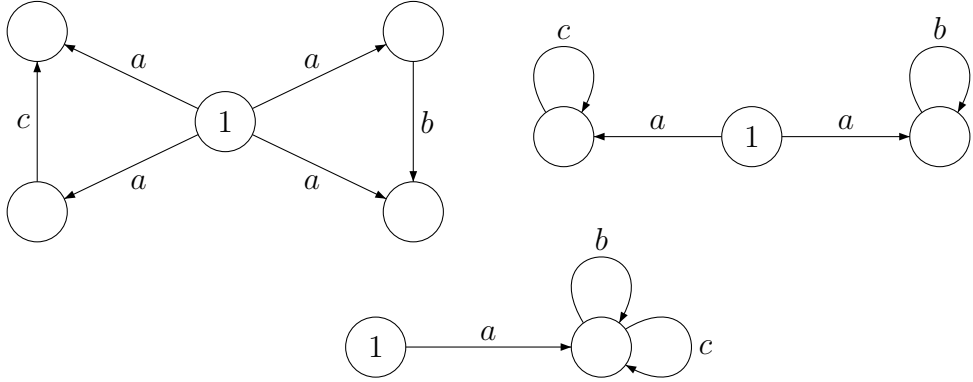
\begin{figure}[htb]
	\centering
    \begin{picture}(120,50)(0,-50)

	\node(n0)(0.0,-2.0){}

    \node(n1)(0.0,-26.0){}

    \node[NLangle=0.0](n2)(25.0,-14.0){1}

    \node(n3)(50.0,-2.0){}

    \node(n4)(50.0,-26.0){}

    \node(n10)(70.0,-16.0){}

    \node[NLangle=0.0](n11)(95.0,-16.0){1}

    \node(n12)(120.0,-16.0){}

    \node[NLangle=0.0](n13)(45.0,-44.0){1}

    \node(n14)(75.0,-44.0){}

    \drawedge[ELside=r](n2,n0){$a$}

    \drawedge(n2,n1){$a$}

    \drawedge(n1,n0){$c$}

    \drawedge(n2,n3){$a$}

    \drawedge[ELside=r](n2,n4){$a$}

    \drawedge(n3,n4){$b$}

    \drawedge[ELside=r](n11,n10){$a$}

    \drawedge(n11,n12){$a$}

    \drawloop(n10){$c$}

    \drawloop(n12){$b$}

    \drawedge(n13,n14){$a$}

    \drawloop(n14){$b$}

    \drawloop[loopangle=0.0](n14){$c$}

    \end{picture}
    \caption{Computing the representation of $H = \langle aba\inv, 
    aca\inv\rangle$}
    \label{fig 1}
\end{figure}

\begin{example}
    Let $A=\{a,b,c\}$.  The above procedure applied to the subgroup
    $H=\langle aba\inv, aca\inv \rangle$ of $F(A)$ is represented in
    Figure~\ref{fig 1}, where the last graph is $\Gamma_A(H)$.
\end{example}
    
Let $\Gamma$ and $\Delta$ be reduced $A$-labeled graphs as above.  A
mapping $\phi$ from the vertex set of $\Gamma$ to the vertex set of
$\Delta$ (we write $\phi\colon\Gamma\rightarrow\Delta$) is a
\textit{morphism of reduced ($A$-)labeled graphs} if it maps the
designated vertex of $\Gamma$ to the designated vertex of $\Delta$ and
if, for each $a\in A$, whenever $\Gamma$ has an $a$-labeled edge $e$
from vertex $u$ to vertex $v$, then $\Delta$ has an $a$-labeled edge
$f$ from vertex $\phi(u)$ to vertex $\phi(v)$.  The edge $f$ is
uniquely defined since $\Delta$ is reduced.  We then extend the domain
and range of $\phi$ to the edge sets of the two graphs, by letting
$\phi(e) = f$.

Note that such a morphism of reduced $A$-labeled graphs is necessarily
locally injective (an \emph{immersion} in \cite{St}), in the following
sense: for each vertex $v$ of $\Gamma$, distinct edges starting (resp.
ending) at $v$ have distinct images.  Further following \cite{St}, we
say that the morphism $\phi\colon\Gamma\rightarrow\Delta$ is a
\emph{cover} if it is locally bijective, that is, if the following
holds: for each vertex $v$ of $\Gamma$, each edge of $\Delta$ starting
(resp.  ending) at $\phi(v)$ is the image under $\phi$ of an edge of
$\Gamma$ starting (resp.  ending) at $v$.

The graph with a single vertex, called 1, and with one $a$-labeled
loop for each $a\in A$ is called the \emph{bouquet of $A$ circles}.
It is a reduced graph, equal to $\Gamma_{A}(F(A))$, and every reduced
graph admits a trivial morphism into it.  One can show that a subgroup
$H$ of $F(A)$ has finite index if and only if this natural morphism
from $\Gamma_{A}(H)$ to the bouquet of $A$ circles is a cover, and in
that case, the index of $H$ in $F(A)$ is the number of vertices of
$\Gamma_A(H)$.  In particular, it is easily decidable whether a
finitely generated subgroup of $F(A)$ has finite index.

This graph-theoretic representation of subgroups of free groups leads
to many more algorithmic results, some of which are discussed at
length in this paper.  We will use some well-known facts
(see~\cite{St}).  If $H$ is a finitely generated subgroup of $F(A)$,
then the rank of $H$ is given by the formula
$$\rank(H) = E - V + 1,$$
where $E$ (resp.  $V$) is the number of edges (resp.  vertices) in
$\Gamma_{A}(H)$.  A more precise result shows how each spanning tree
in $\Gamma_{A}(H)$ (a subtree of the graph $\Gamma_{A}(H)$ which
contains every vertex) determines a basis of $H$.  It is also
interesting to note that if $H$ and $K$ are finitely generated
subgroups of $F(A)$, then $\Gamma_{A}(H\cap K)$ can be easily
constructed from $\Gamma_{A}(H)$ and $\Gamma_{A}(K)$: one first
considers the $A$-labeled graph whose vertices are pairs $(u,v)$
consisting of a vertex $u$ of $\Gamma_{A}(H)$ and a vertex $v$ of
$\Gamma_{A}(K)$, with an $a$-labeled edge from $(u,v)$ to $(u',v')$ if
and only if there are $a$-labeled edges from $u$ to $u'$ in
$\Gamma_{A}(H)$ and from $v$ to $v'$ in $\Gamma_{A}(K)$.  Finally, one
considers the connected component of vertex $(1,1)$ in this product,
and we repeatedly remove the vertices of valence 1, other than the
distinguished vertex $(1,1)$ itself, to make it a reduced $A$-labeled
graph.

To conclude this section, it is very important to observe that if we
change the ambient basis of $F$ from $A$ to $B$, we may radically
modify the labeled graph associated with a subgroup $H$ of $F$, see
Example~\ref{example labeled graphs} below.  In fact, a clearer
understanding of the transformation from $\Gamma_{A}(H)$ to
$\Gamma_{B}(H)$ (put otherwise: of the action of the automorphism
group of $F(A)$ on the $A$-labeled reduced graphs) is one of the
challenges of the field.

\begin{example}\label{example labeled graphs}
    Let $F$ be the free group with basis $A = \{a,b,c\}$, and let $H =
    \langle ab,acba\rangle$.  Note that $B=\{a',b',c'\}$ is also a
    basis of $F$, where $a'=a$, $b'=ab$ and $c'=acba$.  The graphs
    $\Gamma_{A}(H)$ and $\Gamma_B(H)$ are depicted in Figure~\ref{fig
    2}.
\end{example}
    
    \begin{figure}[htb]
	\centering
    \begin{picture}(80,28)(0,-28)

	\node[NLangle=0.0](n0)(0.0,-4.0){1}

    \node(n1)(0.0,-24.0){}

    \node(n2)(30.0,-4.0){}

    \node(n3)(30.0,-24.0){}

    \node[NLangle=0.0](n4)(70.0,-16.0){1}

    \drawedge[curvedepth=3](n0,n2){$a$}

    \drawedge[curvedepth=3](n2,n0){$b$}

    \drawedge(n1,n0){$a$}

    \drawedge(n3,n1){$b$}

    \drawedge(n2,n3){$c$}

    \drawloop[loopangle=20.0,loopdiam=10.0](n4){$c'$}

    \drawloop[loopangle=160.0,loopdiam=10.0](n4){$b'$}

    \end{picture}
    \caption{The graphs $\Gamma_{A}(H)$ and $\Gamma_B (H)$}
    \label{fig 2}
\end{figure}
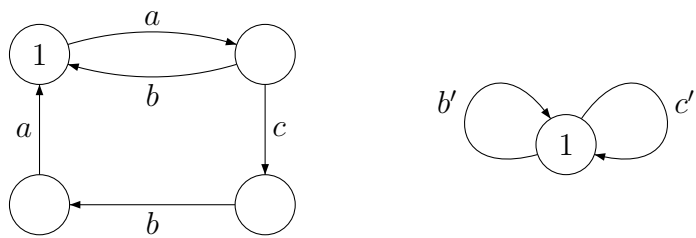

\subsection{Subgroups of subgroups}\label{sub of sub}

A pair of free groups $H\leq K$ is called an \emph{extension} of free
groups.  If $H\leq M\leq K$ are free groups then $H\leq M$ will be
referred to as a \emph{sub-extension} of $H\leq K$.

If $H\le K$ is an extension of free groups, we use the following
shorthand notation: $H\lefg K$ means that $H$ is finitely generated;
$H\lefi K$ means that $H$ has finite index in $K$; and $H\leff K$
means that $H$ is a free factor of $K$.

Extensions can be characterized by means of the labeled graphs
associated with subgroups as in Section~\ref{geometrical apparatus}.
We first note the following simple result (see~\cite[Proposition
2.4]{MSW} or~\cite[Section 4]{KM}).

\begin{lem}
Let $H,K$ be subgroups of a free group $F$ with basis $A$.  Then $H\le
K$ if and only if there exists a morphism of labeled graphs
$\phi_{H,K}$ from $\Gamma_{A}(H)$ to $\Gamma_{A}(K)$.  If it exists,
this morphism is unique.
\end{lem}

Given an extension $H\le K$ between subgroups of the free group with
basis $A$, certain properties of the resulting morphism $\phi_{H,K}$
have a natural translation on the relation between $H$ and $K$.  For
instance, it is not difficult to verify that $\phi_{H,K}$ is a
covering if and only if $H$ has finite index in $K$ (and the index is
the cardinality of each fibre).  This generalizes the characterization
of finite index subgroups of $F(A)$ given in the previous section.

If $\phi_{H,K}$ is one-to-one (and that is, if and only if it is
one-to-one on vertices), then $H$ is a free factor of $K$.
Unfortunately, the converse is far from holding since each non-cyclic
free group has infinitely many free factors.  Furthermore, given $K$,
the particular collection of free factors $H\leff K$ such that
$\phi_{H,K}$ is one-to-one heavily depends on the ambient basis.

We recall here, for further reference, the following well-known
properties of free factors (see~\cite{LS} or~\cite{MKS}).

\begin{lem}\label{standard ppty ff}
    Let $H$, $K$, $L$, $(H_{i})_{i\in I}$ and $(K_{i})_{i\in I}$ be
    subgroups of a free group $F$.
    \begin{itemize}
	\item[(i)] If $H\leff K\leff L$, then $H\leff L$.
	\item[(ii)] If $H_{i} \leff K_{i}$ for each $i\in I$, then
	$\bigcap_{i}H_{i} \leff \bigcap_{i}K_{i}$.
    \end{itemize}
    In particular, if $H$ is a free factor of each $K_{i}$, then $H$
    is a free factor of their intersection; and an intersection of
    free factors of $K$ is again a free factor of $K$.
\end{lem}

Finally, in the situation $H\leq K$, we say that $K$ is an
\emph{$A$-principal overgroup} of $H$ if $\phi_{H,K}$ is onto (both on
vertices and on edges).  We refer to the set of all $A$-principal
overgroups of $H$ as the \textit{$A$-fringe} of $H$, denoted
$\calO_{A}(H)$.  As seen later, this set strongly depends on $A$.  The
$A$-fringe of $H$ is finite whenever $H$ is finitely generated.

Principal overgroups were first considered under the name of
\textit{overgroups} in~\cite{MSW} (see~\cite{Weil00} as well).  They
also appeared later as \emph{principal quotients} in~\cite{KM}, and
their first introduction is in the earlier~\cite{V}, where
$\calO_{A}(H)$ was called the \emph{fringe} of $H$, its \textit{orla}
in catalan.  We shall use the phrase \textit{principal overgroup} (to
stress the fact that not every $K$ containing $H$ is a principal
overgroup of $H$) and \textit{fringe}, omitting the reference to the
basis $A$ when there is no risk of confusion.  Both \textit{orla} and
\textit{overgroup} justify the notation $\calO_{A}(H)$.

Given a finitely generated subgroup $H\leq F(A)$, the fringe
$\calO_{A}(H)$ is computable: it suffices to compute $\Gamma_{A}(H)$,
and to consider each equivalence relation $\sim$ on the set of
vertices of $\Gamma_{A}(H)$.  Say that such an equivalence relation
$\sim$ is a congruence (with respect to the labeled graph structure of
$\Gamma_{A}(H)$) if, whenever $p\sim q$ and there are $a$-labeled
edges from $p$ to $p'$ and from $q$ to $q'$ (resp.  from $p'$ to $p$
and from $q'$ to $q$), then $p'\sim q'$.  Then each congruence gives
rise to a surjective morphism from $\Gamma_{A}(H)$ onto
$\Gamma_{A}(H)/\!\sim$, and hence to a principal overgroup $K$ of $H$
such that $\Gamma_{A}(K) = \Gamma_{A}(H)/\!\sim$.  Moreover, each
principal overgroup $K\in\calO_{A}(H)$ arises in this fashion.  At the
time of writing, a computer program is being developed with the
purpose, among others, of efficiently computing the fringe of a
finitely generated subgroup of a free group (see~\cite{RVW}).

\begin{figure}[htb]
    \centering
\begin{picture}(110,91)(0,-91)

\node(n0)(0.0,-4.0){1}

\node(n1)(0.0,-24.0){4}

\node(n2)(30.0,-4.0){2}

\node(n3)(30.0,-24.0){3}

\drawedge[curvedepth=3](n0,n2){$a$}

\drawedge[curvedepth=3](n2,n0){$b$}

\drawedge(n1,n0){$a$}

\drawedge(n3,n1){$b$}

\drawedge(n2,n3){$c$}

\node(n4)(50.0,-14.0){4}

\node(n5)(80.0,-14.0){1,3}

\node(n6)(110.0,-14.0){2}

\drawedge[curvedepth=3](n4,n5){$a$}

\drawedge[curvedepth=3](n5,n4){$b$}

\drawedge[curvedepth=3](n5,n6){$a$}

\drawedge[curvedepth=2](n6,n5){$b$}

\drawedge[curvedepth=8](n6,n5){$c$}

\node(n7)(0.0,-44.0){1}

\node(n8)(30.0,-44.0){2,4}

\node(n9)(60.0,-44.0){3}

\drawedge[curvedepth=3](n7,n8){$a$}

\drawedge[curvedepth=2](n8,n7){$b$}

\drawedge[curvedepth=8](n8,n7){$a$}

\drawedge[curvedepth=3](n8,n9){$c$}

\drawedge[curvedepth=3](n9,n8){$b$}

\node(n10)(80.0,-44.0){1,3}

\node(n11)(110.0,-44.0){2,4}

\drawedge[curvedepth=2](n10,n11){$a$}

\drawedge[curvedepth=8](n10,n11){$b$}

\drawedge[curvedepth=2](n11,n10){$a$}

\drawedge[curvedepth=8](n11,n10){$b$}

\drawedge[curvedepth=14](n11,n10){$c$}

\node(n12)(10.0,-68.0){1}

\node(n13)(40.0,-68.0){2}

\node(n14)(25.0,-88.0){3,4}

\drawedge[curvedepth=3](n12,n13){$a$}

\drawedge[curvedepth=3](n13,n12){$b$}

\drawedge(n13,n14){$c$}

\drawedge(n14,n12){$a$}

\drawloop[loopangle=0.0](n14){$b$}

\node[NLangle=0.0,Nw=10.0,Nmr=2.0](n15)(80.0,-80.0){{$\scriptstyle 1,2,3,4$}}

\drawloop[loopangle=180.0](n15){$a$}

\drawloop[loopangle=90.0](n15){$b$}

\drawloop[loopangle=0.0](n15){$c$}

\end{picture}
\caption{The six quotients of $\Gamma_{A}(\langle ab,acba\rangle)$}
\label{fig 3}
\end{figure}
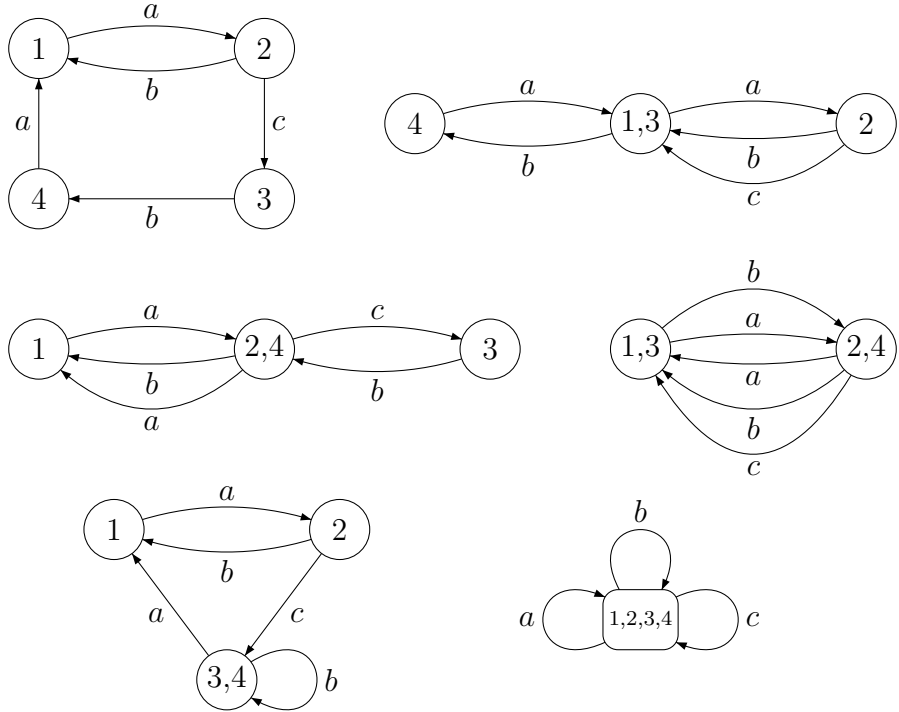

\begin{example}\label{example overgoups}
Let $F$ be the free group with basis $A = \{a,b,c\}$, and let
$H=\langle ab,acba\rangle \leq F$ (the graph $\Gamma_A(H)$ was
constructed in Example~\ref{example labeled graphs}).  Successively
identifying pairs of vertices of $\Gamma_A(H)$ and reducing the
resulting $A$-labeled graph in all possible ways, one concludes that
$\Gamma_A(H)$ has six congruences, whose corresponding quotient graphs
are depicted in Figure~\ref{fig 3}.

Thus the $A$-fringe of $H$ consists on $\calO_{A}(H)=\{H_0, H_1,
H_2, H_3, H_4, H_5\}$, where $H_0=H$, $H_1 =\langle ab, ac,
ba\rangle$, $H_2 =\langle ba, ba\inv, cb\rangle$, $H_3 =\langle ab,
ac, ab\inv, a^2\rangle$, $H_4 =\langle ab, aca, acba\rangle$ and $H_5
=\langle a,b,c\rangle =F(A)$.

However, with respect to the basis $B=\{a, ab, acba\}$ of $F$, the
graph $\Gamma_{B}(H)$ has a single vertex, and hence the $B$-fringe of
$H$ is much simpler, $\calO_B(H) =\{ H\}$.
\end{example}

Finally we observe that, if $H\lefi F(A)$, then $\calO_{A}(H)$
consists of all the extensions of $H$.  Indeed, suppose that $H\le
K\le F(A)$ and $H\lefi F(A)$.  Since $\Gamma_{A}(H)$ is a cover of the
bouquet of $A$ circles (that is, each vertex of $\Gamma_{A}(H)$ is the
origin and the end of an $a$-labeled edge for each $a\in A$), the
range of $\phi_{H,K}$ is also a cover of the bouquet of $A$ circles.
It follows that $\phi_{H,K}$ is onto, since $\Gamma_{A}(K)$ is
connected, and so $K\in \calO_A(H)$.  In particular, if $H\lefi F(A)$,
then $\calO_{A}(H)$ does not depend on $A$, in contrast with what
happens in general.

\subsection{Takahasi's theorem}\label{Takahasi section}

Of particular interest to our discussion is the following 1951 result
by Takahasi (see \cite[Section 2.4, Exercise 8]{MKS}, \cite[Theorem
2]{T} or~\cite[Theorem 1.7]{V}).

\begin{thm}[Takahasi]\label{tak}
Let $F(A)$ be the free group on $A$ and $H\leq F(A)$ a finitely
generated subgroup.  Then, there exists a finite computable collection
of extensions of $H$, say $H=H_0, H_{1},\ldots, H_{n} \leq F(A)$ such
that every extension $K$ of $H$, $H \le K \le F(A)$, is a free
multiple of one of the $H_{i}$.
\end{thm}

The original proof, due to M. Takahasi was combinatorial, using words
and their lengths with respect to different sets of generators.  The
geometrical apparatus described in this section leads to a clear,
concise and natural proof, which was discovered independently by
Ventura in~\cite{V} and by Kapovich and Miasnikov in~\cite{KM}.
Margolis, Sapir and Weil, also independently considered the same
construction in~\cite{MSW} for a slightly different purpose.  Finally,
we note that Turner considered a similar construction in the case of
cyclic subgroups, in his work about test words~\cite{Tu}.  We now give
this proof of Takahasi's theorem.

\demo Let $K$ be an extension of $H$, and let $\phi_{H,K}\colon
\Gamma_{A}(H) \rightarrow \Gamma_{A}(K)$ be the resulting graph
morphism.  Note that the image of $\phi_{H,K}$ is a reduced subgraph
of $\Gamma_A(K)$, and let $L_{H,K}$ be the subgroup of $F(A)$ such
that $\Gamma_{A}(L_{H,K}) = \phi_{H,K}(\Gamma_{A}(H))$.  By
definition, $L_{H,K}$ is an $A$-principal overgroup of $H$ and, by
construction, $\Gamma_{A}(L_{H,K})$ is a subgraph of $\Gamma_A(K)$,
which implies $L_{H,K}\leff K$ (see Section~\ref{sub of sub}).  It
follows immediately that the $A$-fringe of $H$, $\calO_A(H)$,
satisfies the required conditions.
\qed

Thus, for a given $H\lefg F(A)$, the $A$-principal overgroups of $H$
form one possible collection of extensions that satisfy the
requirements of Takahasi's theorem, let us say, a \textit{Takahasi
family for $H$}.  This is certainly not the only one: firstly, we may
add arbitrary subgroups to a Takahasi family; secondly, we observe
that the statement of the theorem does not depend on the ambient
basis, so if $B$ is another basis of $F(A)$, then $\calO_{B}(H)$ forms
a Takahasi family for $H$ as well.  There does however exist a minimum
Takahasi family for $H$ (see Proposition~\ref{AE minimal} below),
which in particular does not depend on the ambient basis.  The main
object of this paper is a discussion of this family, which is
introduced in the next section.

\section{Algebraic extensions}\label{sec algebraic extensions}

The notion of algebraic extension discussed in this paper was first
introduced by Kapovich and Miasnikov \cite{KM}.  It seems to be mostly
of interest for finitely generated subgroups, but many definitions and
results hold in general and we avoid restricting ourselves to finitely
generated subgroups until that becomes necessary.

\subsection{Definitions}\label{defs section}

Let $H\leq K$ be an extension of free groups and let $x\in K$.  We say
that $x$ is \textit{$K$-algebraic over $H$} if every free factor of
$K$ containing $H$, $H\leq L\leff K$, satisfies $x\in L$.  Otherwise
(i.e. if there exists $H\leq L\leff K$ such that $x\not\in L$) we say
that $x$ is \textit{$K$-transcendental over $H$}.

\begin{example}\label{example algebraic elements}
If $H\le K$, then every element $x\in H$ is obviously $K$-algebraic
over $H$.

Every element $x\in K$ is $K$-algebraic over $\langle x^n\rangle$, for
each integer $n\neq 0$.  In fact, it is straightforward to verify that
if $x^n$ lies in a free factor $L$ of $K$, then so does $x$.

If $x$ is primitive in $K$ (that is, if $\langle x \rangle \leff K$),
then every element of $K \setminus \langle x\rangle$ is
$K$-transcendental over the subgroup $\langle x\rangle$.

The notion of algebraicity over $H$ is relative to $K$.  For example,
in $F = F(a,b)$, $a^2$ is $\langle a^2, b^2 \rangle$-transcendental
over $H = \langle a^2 b^2 \rangle$ since $a^2b^2$ is primitive in
$\langle a^2, b^2 \rangle$.  However, $a^2$ is $F$-algebraic over $H$
because no proper free factor of $F$ contains $a^2b^2$.
\end{example}

The following is a trivial but useful observation.

\begin{obs}\label{product}
     Let $H\leq K$ be an extension of free groups, and let $x,y\in K$.
     \begin{itemize}
        \item[(i)] If $x,y$ are $K$-algebraic over $H$ then so are
        $x^{-1}$ and $xy$.

        \item[(ii)] If $x,y$ are $K$-transcendental over $H$ then so is
        $x^{-1}$ (but not in general $xy$).
     \end{itemize}
\end{obs}

We say that an extension of free groups $H\leq K$ is \emph{algebraic},
and we write $H\lealg K$, if every element of $K$ is $K$-algebraic
over $H$.  It is called \emph{purely transcendental} if every element
of $K$ is either in $H$ or is $K$-transcendental over $H$.  Naturally,
there are extensions that are neither algebraic nor purely
transcendental.  These concepts were originally introduced
in~\cite{KM}, and the following propositions further describe their
properties.

\begin{prop}\label{algext}
Let $H\leq K$ be an extension of free groups.  The following are
equivalent:
\begin{itemize}
\item[(a)] $H$ is contained in no proper free factor of $K$;
\item[(b)] $H\lealg K$, that is, every $x\in K$ is $K$-algebraic over
$H$;
\item[(c)] there exists $X\subseteq K$ such that $K = \langle H\cup
X\rangle$ and every $x\in X$ is $K$-algebraic over $H$ (furthermore,
if $K$ is finitely generated, one may choose $X$ to be finite).
\end{itemize}
\end{prop}

\demo (b) follows from (a) by definition.  If (b) holds, then (c)
holds with $X$ any system of generators for $K$.  Finally, (a) follows
from (c) in view of Fact~\ref{product}~(i).  \qed

\begin{prop}\label{ptrext}
Let $H\leq K$ be an extension of free groups. The following are equivalent:
\begin{itemize}
\item[(a)] $H$ is a free factor of $K$,
\item[(b)] $H\le K$ is purely transcendental, that is, every $x\in
K\setminus H$ is $K$-transcendental over $H$.
\end{itemize}
\end{prop}

\demo (a) implies (b) by definition.  To prove the converse, let $M$
be the intersection of all the free factors of $K$ containing $H$.  By
Lemma~\ref{standard ppty ff}, $M$ is a free factor of $K$ containing
$H$, and (b) implies that $M=H$.  \qed

\begin{example}\label{example alg ext}
It is easily verified (say, using Example~\ref{example algebraic
elements}) that if $1\neq x\in F$ and $n\neq 0$, we have $\langle
x^n\rangle \lealg \langle x\rangle$.

By Proposition~\ref{ptrext}, an extension of the form $\langle
x\rangle \leq F$ is purely transcendental if and only if $x$ is a
primitive element of $F$.  Moreover, if $F$ has rank two, then
$\langle x\rangle \leq F$ is algebraic if and only if $x$ is not a
power of a primitive element of $F$.

Assuming again that $F$ has rank two, $H\lealg F$ for every non-cyclic
subgroup $H$.  Indeed, every proper free factor of $F$ is cyclic and
hence cannot contain $H$.
\end{example}

We denote by $\AE(H)$ the set of algebraic extensions of $H$, and we
observe that, in contrast with the definition of principal overgroups,
this set does not dependent on the choice of an ambient basis.  This
same observation can be expressed as follows.

\begin{obs}
Let $H\le K\le F$ be extensions of free groups and let
$\phi\in\Aut(F)$.  Then $H\lealg K$ if and only if $\phi(H) \lealg
\phi(K)$.
\end{obs}

We can now express the connection between algebraic extensions and
Takahasi's theorem.

\begin{prop}\label{AE minimal}
Let $H\lefg F(A)$ be an extension of free groups. Then we have:
\begin{enumerate}
\item[(i)] $\AE(H)\subseteq \calO_{A}(H)$;
\item[(ii)] $\AE(H)$ is finite (i.e., $H$ admits only a finite number
of algebraic extensions);
\item[(iii)] $\AE(H)$ is the set of $\leff$-minimal elements of every
Takahasi family for $H$ (see Section~\ref{Takahasi section});
\item[(iv)] $\AE(H)$ is the minimum Takahasi family for $H$.
\end{enumerate}
\end{prop}

\demo Let $K$ be an algebraic extension of $H$.  The proof of
Takahasi's theorem shows that $K$ is a free multiple of some principal
overgroup $L\in \calO_A(H)$.  Then, Proposition~\ref{algext} implies
that $L=K$ proving (i).  Statement (ii) follows immediately.

Let $\calL$ be a Takahasi family for $H$ and let $K\in \AE(H)$.  By
definition of $\calL$, there exists a subgroup $L\in\calL$ such that
$H\le L \leff K$.  By Proposition~\ref{algext}, it follows that $L =
K$, so $K\in\calL$.  Thus $\AE(H)$ is contained in every Takahasi
family for $H$.  For the same reason, $K$ is $\leff$-minimal in
$\calL$.

Now suppose that $K\in \calL$ is $\leff$-minimal in $\calL$, and let
$M$ be an extension of $H$ such that $H\le M\leff K$.  By definition
of a Takahasi family, there exists $L\in\calL$ such that $H\le L\leff
M$, so $L\leff K$.  Since $K$ is $\leff$-minimal in $\calL$, it
follows that $L=K$, so $M=K$.  Hence, $H\lealg K$ concluding the proof
of (iii).

Finally, it is immediate that the $\leff$-minimal elements of a
Takahasi family for $H$ again form a Takahasi family.  Statement (iv)
follows directly.  \qed

\begin{example}\label{alg ext vs fi}
If $H\lefi K$, then $H\lealg K$.  This follows immediately from the
observation that a proper free factor of $K$ has infinite index.

It follows that, if $H\lefi F(A)$, then $\AE(H) = \calO_{A}(H)$ is
equal to the set of all extensions of $H$.  Indeed, we have already
observed at the end of Section~\ref{sub of sub} that every extension
of $H$ is an $A$-principal overgroup of $H$, and since $H$ has finite
index in each of its extensions, it is algebraic in each.
\end{example}

Proposition~\ref{AE minimal} shows that $\AE(H)$ is contained in
$\calO_{A}(H)$ for each ambient basis $A$.  We conjecture that
$\AE(H)$ is in fact equal to the intersection of the sets
$\calO_{A}(H)$, when $A$ runs over all the bases of $F$.
Example~\ref{alg ext vs fi} shows that the conjecture holds if $H$ has
finite index.  It also holds if $H\leff F$, since in that case,
$\AE(H) = \{H\}$, and $F$ admits a basis $B$ relative to which
$\Gamma_{B}(H)$ is a graph with a single vertex.

We conclude with a simple but important statement.

\begin{prop}\label{computability AE}
Let $F(A)$ be the free group on $A$ and $H\lefg F(A)$.  The set
$\AE(H)$ is computable.
\end{prop}

\demo Since every algebraic extension of $H$ is in $\calO_{A}(H)$, it
suffices to compute $\calO_{A}(H)$ and then, for each pair of distinct
elements $K,L\in\calO_{A}(H)$, to decide whether $L\leff K$: $\AE(H)$
consists of the principal overgroups of $H$ that do not contain
another principal overgroup as a free factor.

In order to conclude, we observe that deciding whether $L\leff K$ can
be done, for example, using the first part the classical Whitehead's
algorithm.  More precisely, Whitehead's algorithm (see
\cite[Proposition 4.25]{LS}) shows how to decide whether a tuple of
elements, say $u = (u_{1},\ldots,u_{r})$, of a free group $K$ can be
mapped to another tuple $v=(v_1, \ldots ,v_r)$ by some automorphism of
$K$.  The first part of this algorithm reduces the sum of the length
of the images of the $u_{i}$ to its minimal possible value.  And it is
easy to verify that this minimal total length is exactly $r$ if and
only if $\{ u_1, \ldots ,u_r\}$ freely generates a free factor of $K$.
We point out here that an alternative algorithm was recently proposed
by Silva and Weil~\cite{SW}.  That algorithm is faster, and completely
based on graphical tools.  \qed

The efficiency of the algorithm to compute $\AE(H)$ sketched in the
proof of Proposition~\ref{computability AE}, is far from optimal.  An
upcoming paper by A. Roig, E. Ventura and P. Weil discusses better
computation techniques for that purpose~\cite{RVW}.

\begin{remark}
The terminology adopted for the concepts developed in this section is
motivated by an analogy with the theory of field extensions.  More
precisely, if an element $x \in K$ is $K$-transcendental over $H$,
then $H$ is a free factor of $\langle H,x\rangle$ and $\langle
H,x\rangle =H*\langle x\rangle$ (see Proposition~\ref{elext} below).
This is similar to the field-theoretic definition of transcendental
elements: an element $x$ is transcendental over $H$ if and only if the
field extension of $H$ generated by $x$ is isomorphic to the field of
rational fractions $H(X)$.

However, the analogy is not perfect and in particular, the converse
does not hold.  For instance, $a^2$ is $\langle a,b\rangle$-algebraic
over $\langle a^2b^2\rangle$ (see Example~\ref{example algebraic
elements}), but $\langle a^2b^2, a^2\rangle = \langle a^2b^2\rangle *
\langle a^2\rangle$.  This stems from the fact, noticed earlier, that
the notion of an element $x$ being $K$-algebraic over $H$, depends on
$K$ and not just on $x$.

It is natural to ask whether the analogy also extends to the
definition of algebraic elements: in other words, is there a natural
analogue in this context for the notion of roots of a polynomial with
coefficients in $H$?  The discussion of equations in Section~\ref{sec
equations} offers some insight into this question.
\end{remark}

\subsection{Composition of extensions}

We now consider compositions of extensions.  Some of the results in
the following proposition come from~\cite{KM}.  We restate and extend
them here with simpler proofs.  We also include in the statement
well-known facts (the primed statements), in order to emphasize the
dual properties of algebraic and purely transcendental extensions.

\begin{prop}\label{composition}
Let $H\leq K$ be an extension of free groups, and let $H\leq K_i \leq
K$ be two sub-extensions, $i=1,2$.
\begin{itemize}
\item[(i)] If $H\lalg K_1\lalg K$ then $H\lalg K$.
\item[(i')] If $H\lff K_1\lff K$ then $H\lff K$.
\item[(ii)] If $H\lalg K$ then $K_1 \lalg K$, while $H\leq K_1$ need
not be algebraic.
\item[(ii')] If $H\lff K$ then $H\lff K_1$, while $K_1 \leq K$ need
not be purely transcendental.
\item[(iii)] If $H\lalg K_1$ and $H\lalg K_2$ then $H\lalg \langle K_1
\cup K_2\rangle$, while $H\leq K_1\cap K_2$ need not be algebraic.
\item[(iii')] If $H\lff K_1$ and $H\lff K_2$ then $H\lff K_1\cap K_2$,
while $H\leq \langle K_1 \cup K_2\rangle$ need not be purely
transcendental.
\end{itemize}
\end{prop}

\demo Statement \textit{(i')} and the positive parts of statements
\textit{(ii')} and \textit{(iii')} can be found in Lemma~\ref{standard
ppty ff}.  The free group $F$ on $\{a,b\}$ already contains
counterexamples for the converse statements in \textit{(ii')} and
\textit{(iii')}: for the first one, we have $\langle a \rangle \lff F$
while $\langle a \rangle \leff \langle a, b^2 \rangle \lalg F$ (see
Example~\ref{example alg ext}).  And for the second one, we have
$\langle[a,b]\rangle \leff \langle a, [a,b]\rangle$ and
$\langle[a,b]\rangle \leff \langle b, [a,b]\rangle$, whereas
$\langle[a,b]\rangle \lealg \langle a, [a,b], b\rangle = F$.

Now assume that $H\lalg K_1\lalg K$ and let $L$ be a free factor of
$K$ containing $H$.  Then, $L\cap K_1$ is a free factor of $K_1$
containing $H$ by Lemma~\ref{standard ppty ff}.  Since $H\leq K_1$ is
algebraic, we deduce that $L\cap K_1 = K_1$, and hence $K_1 \leq L$.
But $K_1 \lealg K$, so $L=K$.  Thus, the extension $H\leq K$ is
algebraic, which proves \textit{(i)}.

The first part of \textit{(ii)} is clear.  A counterexample for the
second part in $F = F(a,b)$ is as follows: we have $\langle [a,b]
\rangle \leff \langle a, [a,b] \rangle \le F$, while $\langle [a,b]
\rangle \lealg F$ by Example~\ref{example alg ext}.

Suppose now that $H\lalg K_1$ and $H\lalg K_2$, and let $L$ be a free
factor of $\langle K_1 \cup K_2 \rangle$ containing $H$.  Then
Lemma~\ref{standard ppty ff} shows that, for $i=1,2$, $L\cap K_i \leff
K_i$ containing $H$.  Since $H\lalg K_i$, we deduce that $L\cap K_i =
K_i$ and hence, $K_i \leq L$.  Thus, $L = \langle K_1 \cup K_2
\rangle$, and the extension $H\leq \langle K_1 \cup K_2 \rangle$ is
algebraic, thus proving the positive part of \textit{(iii)}.

Finally, to conclude the proof of \textit{(iii)}, it suffices to
exhibit subgroups $H$, $K_{1}$, $K_{2}$ such that $H\lealg K_{i}$
($i=1,2$) but $H\leff K_{1}\cap K_{2}$.  Again in $F(a,b)$ take, for
example, $K_{1}=\langle a^2,b \rangle$ and $K_{2}=\langle a^3,
b\rangle$, whose intersection is $K_1 \cap K_2 =\langle a^6,
b\rangle$.  Letting $H=\langle a^6b \rangle$, we have $H\leff
K_{1}\cap K_{2}$ but $H\lealg K_1$ and $H\lealg K_2$.  \qed

To close this section, let us note another natural property of
algebraic extensions, which slightly generalizes a result of Kapovich
and Miasnikov \cite{KM}.

\begin{prop}\label{join of algebraic}
Let $F$ be a free group.  If $H_i \lalg K_i \leq F$ ($i\in I$), then
$\langle \bigcup_{i}H_i \rangle \lalg \langle \bigcup_{i}K_i \rangle$.
The converse holds if $\langle \bigcup_{i}K_i\rangle = *_{i} K_{i}$.
\end{prop}

\demo Suppose that $\langle \bigcup_{i}H_i \rangle \le L \lff \langle
\bigcup_{i}K_i \rangle$.  Let $j\in I$.  By Lemma~\ref{standard ppty
ff}, we have $L\cap K_j \lff \langle \bigcup_{i}K_i \rangle \cap K_j =
K_j$.  Moreover, $H_{j}\le L\cap K_{j}$, so $K_{j} = L\cap K_{j}$
since $H_{j}\lealg K_{j}$, and hence $K_{j}\subseteq L$.  This holds
for each $j\in I$, so $L = \langle \bigcup_{i}K_i \rangle$ and we have
shown that $\langle \bigcup_{i}H_i \rangle \lalg \langle
\bigcup_{i}K_i \rangle$.

For the converse, suppose that $\langle \bigcup_{i}K_i\rangle =
*_{i}K_{i}$.  It follows that $\langle \bigcup_{i}H_i\rangle =
*_{i}H_{i}$.  Now we assume that $*_i H_i \lealg *_{i} K_i$.  Let
$j\in I$.  If $H_{j}\le L\leff K_{j}$, then $*_{i}H_{i} \le L *
*_{i\ne j}K_{i} \leff *_{i}K_{i}$ and hence $L * *_{i\ne j}K_{i} =
*_{i}K_{i}$.  Taking the projection onto $K_{j}$, it follows that $L =
K_{j}$.  Thus $H_{j}\lealg K_{j}$ for each $j\in I$.  \qed

Note that the converse of Proposition~\ref{join of algebraic} does not
hold in general, as can be seen from the counterexample provided in
the proof of Proposition~\ref{composition} $(iii')$.

\subsection{Elementary extensions}

We say that an extension of free groups $H\leq K$ is \emph{elementary}
if $K=\langle H, x \rangle$ for some $x\in K$.  Elementary extensions
turn out to be either algebraic or purely transcendental, as we now
see.

\begin{prop}\label{elext}
Let $H \le F$ be an extension of free groups and let $x\in F$.  Let
also $X$ be a new letter, not in $F$.  The following are equivalent:
\begin{itemize}
\item[(a)] the morphism $H * \langle X \rangle \to F$ acting as the
identity over $H$ and sending $X$ to $x$ is injective;
\item[(b)] $H$ is a proper free factor of $\langle H, x \rangle$;
\item[(c)] $H$ is contained in a proper free factor of $\langle H, x
\rangle$.
\end{itemize}
If, in addition, $H$ is finitely generated, then these are further
equivalent to:
\begin{itemize}
\item[(d)] $\rank(\langle H,x \rangle)=\rank(H)+1$;
\item[(e)] $\rank(\langle H,x \rangle)>\rank(H)$.
\end{itemize}
\end{prop}

\demo It is immediately clear that statement~(a) implies~(b), and
that~(b) implies~(c).

At this point, let us assume that $H$ has finite rank.  It is
immediate that $\rank(\langle H,x\rangle)\le \rank(H)+1$, so~(b)
implies~(d) and~(d) and~(e) are equivalent.  Now consider the morphism
from $H*\langle X\rangle$ to $\langle H, x\rangle$ mapping $H$
identically to itself, and $X$ to $x$.  This morphism is surjective by
construction, and if $\rank(\langle H, x\rangle) = \rank(H*\langle
X\rangle)$, then it is injective by the hopfian property of finitely
generated free groups.  That is, (d) implies~(a).  Thus we have shown
that if $H$ has finite rank, then statements~(a), (b), (d) and~(e) are
equivalent.  It only remains to prove that~(c) implies~(a).

We now return to the general case, where $H$ may have infinite rank,
and we assume that~(c) holds, that is, $\langle H, x\rangle = K*L$ for
some $L\neq 1$ and $H\leq K$.  We have $\langle H,x\rangle \le \langle
K, x\rangle$, and hence $\langle K, x\rangle = \langle H, x\rangle =
K*L$.  Moreover, $x\not\in K$ and we let $x=k_0 \ell_1 k_1\cdots
\ell_r k_r$ be the normal form of $x$ in the free product $K*L$.

Let $M$ be a finitely generated free factor of $K$ containing the
$k_{i}$, and let $N$ be such that $K=M*N$.  First we observe that
$$
\langle M, x\rangle \le M * L \leff K * L = M * N * L = \langle K,
x\rangle = \langle M, N, x\rangle.
$$
It follows that $N$ is a free complement of $\langle M, x\rangle$ in
$\langle K, x\rangle$, that is, $\langle K, x\rangle =\langle M,
x\rangle *N$.

Next we note that $M\leff K\leff \langle K, x\rangle$, so $M \leff
\langle K, x\rangle$ and hence $M \leff \langle M, x\rangle$.  Since
$x\not\in K$, $M$ is a finitely generated, proper free factor of
$\langle M, x\rangle$, and we already know that this implies that the
morphism from $M*\langle X\rangle$ to $F$ mapping $M$ identically to
itself and mapping $X$ to $x$, is injective.  Since $N$ is a free
complement of the range of this morphism in $\langle H, x\rangle$, and
also a free complement of $M$ in $K$, it follows that the natural
mapping from $K*\langle X\rangle$ to $F$ mapping $X$ to $x$ is
injective.  Its restriction to $H*\langle X\rangle$ is therefore
injective, and statement~(a) holds, which completes the proof.  \qed

Proposition~\ref{elext} immediately translates into the following.

\begin{cor}\label{+1}
Let $F$ be a free group and $H \le K$ be an elementary extension of
subgroups of $F$.  Then, either $H\lalg K$ or $H\lff K$.  Furthermore,
if $H$ is finitely generated then $\rank(K)\leq \rank(H)+1$ with
equality if and only if $H\lff K$.
\end{cor}

Let us say that an extension $H\le K$ is \emph{e-algebraic}, written
$H\leealg K$, if it splits as a finite composition of algebraic,
elementary extensions, $H\lalg H_1 \lalg \cdots \lalg H_k = K$.  Then
Proposition~\ref{elext} yields the following.

\begin{cor}\label{rk goes down if elem*}
Let $H$ be a finitely generated subgroup of a free group $F$ and let
$H\leealg K$ be an e-algebraic extension.  Then $\rank(K)\leq
\rank(H)$.
\end{cor}

Obviously, every extension $H\leq K$ with $K$ finitely generated,
splits into a composition of elementary extensions, but an algebraic
extension $H\lealg K$ cannot always be split into a composition of
algebraic elementary extensions.  In view of Corollary~\ref{rk goes
down if elem*}, this is the case for the algebraic extension $\langle
[a,b] \rangle \lealg F(a,b)$.  Thus, $H\lealg K$ does not imply
$H\leealg K$.

\subsection{Algebraic closure of a subgroup}\label{sec alg closure}

If $H\le K$ is an extension of free groups, there exists a greatest
algebraic extension of $H$ inside $K$.  This can be deduced from
Proposition~\ref{join of algebraic}, but the following theorem is a
more precise statement.

\begin{thm}\label{existcl}
Let $H\leq L\leq K$ be extensions of free groups.  The following are
equivalent.
\begin{itemize}
\item[(a)] $H\lealg L\leff K$.  \item[(b)] $L$ is the intersection of
the free factors of $K$ containing $H$.  \item[(c)] $L$ is the set of
elements of $K$ that are $K$-algebraic over $H$.  \item[(d)] $L$ is
the greatest algebraic extension of $H$ contained in $K$.
\end{itemize}
In this case, the subgroup $L$ is uniquely determined by $H$ and $K$.
\end{thm}

\demo Let $x\in K$.  By definition, $x$ is $K$-algebraic over $H$ if
and only if $x$ sits in every free factor of $K$ containing $H$.  This
is exactly the equivalence of statements~(b) and~(c).  The equivalence
of~(c) and~(d) is a direct consequence of the fact that the elements
that are $K$-algebraic over $H$ form a subgroup (Fact~\ref{product}).
Thus statements~(b), (c) and~(d) are equivalent.

Now let $L$ be defined as in~(b): by~(d), $H\lealg L$.  Now let $x\in
K \setminus L$.  Since $x$ is not algebraic over $H$, there exists a
free factor $M\leff K$ containing $H$ and missing $x$.  But $L\le M$,
so $x$ is not $K$-algebraic over $L$ either.  It follows that the
extension $L\le K$ is purely transcendental, and hence $L\leff K$ by
Proposition~\ref{ptrext}.  This proves~(b) implies~(a).

Finally, let us assume that $H\lealg L \leff K$ for some $L$.  Let $M$
be such that $H\le M\leff K$.  Then $L\cap M \leff L$ by
Lemma~\ref{standard ppty ff}~(ii).  But we also have $H\le L\cap M \le
L$ and $H\lealg L$.  It follows that $L\cap M = L$, that is $L\leq M$,
and~(b) follows.  This concludes the proof.  \qed

\begin{remark}
    It is interesting to compare Theorem~\ref{existcl} with M. Hall's
    Theorem, stating that every finitely generated subgroup $H\leq F$
    is a free factor of a subgroup $M$ of finite index in $F$.  In
    other words, one can split the extension $H\leq F$ in two parts,
    $H\lff M\lfi F$, the first being purely transcendental, and the
    second being finite index (and hence, algebraic).  Note that the
    intermediate subgroup $M$ is not unique in general.
    Theorem~\ref{existcl} yields a ``dual" splitting of the extension
    $H\le F$, where the order between the transcendental and the
    algebraic parts is switched around, and with the additional nice
    property that the intermediate extension is now uniquely
    determined by $H\le F$.
\end{remark}

Let $H\le K$ be an extension of free groups.  The subgroup $L$
characterized in Theorem~\ref{existcl} is called the
\emph{$K$-algebraic closure} of $H$, denoted $\cl_K (H)$.  It is
natural to consider the extremal situations, where $\cl_K(H) = H$ (we
say that $H$ is $K$-\emph{algebraically closed}) and where $\cl_K(H) =
K$ (we say that $H$ is $K$-\emph{algebraically dense}).  Of course,
these situations coincide with $H\leq K$ being purely transcendental
and algebraic, respectively.

\begin{obs}
Let $H\leq K$ be an extension of free groups. Then,
\begin{itemize}
\item[(i)] $H$ is $K$-algebraically closed  if and only if $H\lff K$,
\item[(ii)] $H$ is $K$-algebraically dense if and only if $H\lalg K$.
\end{itemize}
\end{obs}

As established in the following proposition, maximal proper retracts
of a finitely generated free group $K$ are good examples of extremal
subgroups, i.e. subgroups of $K$ that are either $K$-algebraically
closed or $K$-algebraically dense.  Recall that a subgroup $H\leq K$
is a \emph{retract} of $K$ if the identity $\id\colon H\to H$ extends
to a homomorphism $K\to H$, called a \textit{retraction}
(see~\cite{MKS} for a general description of retracts of finitely
generated free groups); in particular, free factors of $K$ are
retracts of $K$.  Note that if $H$ is a retract of $K$ then
$\rank(H)\leq \rank(K)$.  Moreover, if $K$ is finitely generated, the
hopfian property of finitely generated free groups shows that $K$ is
the unique retract of $K$ with rank equal to $\rank(K)$.  So, if $H$
is a proper retract of $K$ then $\rank(H)<\rank(K)$.

We also say that $H$ is \emph{compressed in $K$} (see~\cite{DV}) if
$\rank(H) \le \rank(L)$ for each $H\le L\le K$.  By restricting a
retraction to $L$, it is clear that every retract of $K$ (and, in
particular, every free factor of $K$) is compressed in $K$.

\begin{prop}
Let $K$ be a finitely generated free group.  A maximal proper
compressed subgroup (resp.  a maximal proper retract) $H$ of $K$ is
either $K$-algebraically dense, or $K$-algebraically closed.  In the
latter case, $H$ is in fact a free factor of $K$, of rank
$\rank(K)-1$.
\end{prop}

\demo The algebraic closure $\cl_K(H)$ is a free factor of $K$, and
hence it is also a retract and a compressed subgroup.  By definition
of $H$, either $\cl_{K}(H) = K$, and $H$ is $K$-algebraically dense;
or $\cl_{K}(H)=H$ and $H$ is $K$-algebraically closed and a free
factor.  Maximality then implies the announced rank property.  \qed

We now discuss the behavior of the algebraic closure operator.

\begin{prop}
Let $H_i\leq K$, $i=1,2$, be two extensions of free groups.  Then,
$\cl_{K}(H_1 \cap H_2)\leff \cl_K(H_1) \cap \cl_K(H_2)$, and the
equality is not true in general.
\end{prop}

\demo By Theorem~\ref{existcl}, $\cl_{K}(H_{i})$ is a free factor of
$K$ containing $H_{i}$, so $\cl_K(H_1) \cap \cl_K(H_2)$ is a free
factor of $K$ containing $H_1 \cap H_2$ (Lemma~\ref{standard ppty
ff}).  Again by~Theorem~\ref{existcl}, $\cl_{K}(H_{1}\cap H_{2})$ is a
free factor of $\cl_K(H_1) \cap \cl_K(H_2)$.

A counterexample for the reverse inclusion is as follows: let $K =
F(a,b)$, $H_1 =\langle [a,b]\rangle$ and $H_2 =\langle
[a,b^{-1}]\rangle$.  Both these subgroups are $K$-algebraically dense
(see Example~\ref{example alg ext}) and their intersection is trivial.
\qed

\begin{prop}\label{int2}
Let $K_{i}\leq K$, $i=1,2$, be two extensions of free groups and let
$H\le K_{1}\cap K_{2}$.  Then, $\cl_{K_1\cap K_2}(H)\leff
\cl_{K_1}(H)\cap \cl_{K_2}(H)$, and the equality is not true in
general.
\end{prop}

\demo By Theorem~\ref{existcl}, $\cl_{K_{i}}(H)$ is a free factor of
$K_{i}$ containing $H$, so $\cl_{K_{1}}(H) \cap \cl_{K_2}(H)$ is a
free factor of $K_{1}\cap K_{2}$ containing $H$ (Lemma~\ref{standard
ppty ff}).  Again by~Theorem~\ref{existcl}, $\cl_{K_{1}\cap K_{2}}(H)$
is a free factor of $\cl_{K_{1}}(H) \cap \cl_{K_2}(H)$.

The following is a counter-example for the converse inclusion.  Let $K
= \langle a,b,c \rangle$ be a free group of rank 3, let $H = \langle
[a,b],[a,c]\rangle$, $K_1 = \langle a,b,[a,c]\rangle$ and $K_2
=\langle a,c,[a,b]\rangle$.  One can verify that $K_{1}\cap K_{2} =
\langle a, [a,b], [a,c]\rangle$, so $\cl_{K_{1}\cap K_{2}}(H) = H$.
On the other hand, $H\lealg K_{i}$ by Example~\ref{example alg ext}
and Proposition~\ref{join of algebraic}, so $\cl_{K_{i}}(H) = K_{i}$
and $\cl_{K_{1}}(H) \cap \cl_{K_2}(H) = K_{1}\cap K_{2} \neq
\cl_{K_{1}\cap K_{2}}(H)$.  \qed

\begin{remark}
If $H\le K_{1} \le K_{2}$, Proposition~\ref{int2} shows that
$\cl_{K_{1}}(H) \le \cl_{K_{2}}(H)$.  If in addition $K_1 \leff K_2$,
Proposition~\ref{existcl} shows that $\cl_{K_{1}}(H) =
\cl_{K_{2}}(H)$.  However, in general, even the inclusion
$\cl_{K_1}(H)\leq K_1 \cap \cl_{K_2}(H)$ may be strict, as the
following counterexample shows.

Let $K_2 =\langle a,b\rangle$ be a free group of rank 2, and let $H =
\langle [a,b]\rangle$ and $K_1 = \langle a, [a,b]\rangle$.  Then
$H\lff K_1\lalg F$ and $H\lalg F$.  So, $\cl_{K_{1}}(H) = H$ is
properly contained in $K_{1}\cap \cl_{F}(H) = K_{1}\cap F =K_1$.
\end{remark}

Finally, let us consider e-algebraic extensions.  There too, there
exists a greatest e-algebraic extension, at least for finitely
generated subgroups.  We first prove the following technical lemma.

\begin{lem}\label{technical e-algebraic}
Let $H \le K\le F$ be extensions of free groups and let $x\in F$.  If
$H\lealg \langle H, x\rangle$, then $K\lealg \langle K,x\rangle$.
\end{lem}

\demo Assume $H\lealg \langle H, x\rangle$.  If $K\le \langle K,
x\rangle$ is not algebraic, then $x\not\in K$ and $K \leff \langle K,
x\rangle$ by Proposition~\ref{elext}.  It follows that $H\leq \langle
H,x\rangle \cap K \leff \langle H,x\rangle \cap \langle K,x\rangle =
\langle H,x\rangle$, which forces either $\langle H,x\rangle \cap K
=H$ or $\langle H,x\rangle \cap K =\langle H, x\rangle$.  The first
possibility implies $H=\langle H,x\rangle \cap K \lff \langle
H,x\rangle$ contradicting the hypothesis, while the second possibility
contradicts $x\not\in K$.  \qed

\begin{cor}\label{exists ealg closure}
Let $H\le F$ be an extension of free groups and let $H\leealg K_{i}$
($i = 1,\ldots,n$) be a finite family of e-algebraic extensions of
$H$.  Then $K_{i}\leealg \langle \bigcup_{j}K_{j}\rangle$ for each
$i$.

In particular, if $H$ is finitely generated, then $H$ admits a
greatest e-algebraic extension in $F$.
\end{cor}

\demo It suffices to prove the first statement for $n= 2$.  Let us
assume that $H = H_{0} \lealg H_{1} \lealg \cdots \lealg H_{p} =
K_{1}$ and that $x_{1},\ldots, x_{p}$ are such that $H_{i} = \langle
H_{i-1}, x_{i}\rangle$ for each $1\le i\le p$.  Then a repeated
application of Lemma~\ref{technical e-algebraic} shows that $K_{2}
\leealg \langle K_{2}, x_{1},\ldots, x_{p}\rangle = \langle K_{1}\cup
K_{2}\rangle$.

If $H\lefg F$, $H$ has finitely many algebraic extensions, and among
them finitely many e-algebraic extensions.  The join of these
extensions is again an e-algebraic extension and this concludes the
proof.  \qed

The greatest e-algebraic extension of a subgroup $H\le F$, whose
existence is asserted in Corollary~\ref{exists ealg closure}, is
called its \emph{e-algebraic closure}.  We say that $H$ is
\textit{e-algebraically closed} if it is equal to its e-algebraic
closure.  Proposition~\ref{elext} immediately implies the following
characterization.

\begin{cor}\label{charact ealg closed}
Let $H\le F$ be an extension of free groups.  Then $H$ is
e-algebraically closed if and only if $\langle H,x\rangle =H*\langle
x\rangle$ for each $x\not\in H$.
\end{cor}

\begin{example}
    Let $x\in F$ be an element of a free group not being a proper
    power.  Then, for every $y\in F$, either $\langle x\rangle
    =\langle x,y\rangle$ or $\rank (\langle x,y \rangle )=2$.  In
    other words, maximal cyclic subgroups of free groups are
    e-algebraically closed.
    
    A subgroup $H\le F$ is said to be \textit{strictly compressed} if
    $\rank(H)<\rank(K)$ for each proper extension $H<K\leq F$.  It is
    immediate that strictly compressed subgroups form a natural class
    of e-algebraically closed subgroups.
    
    By Example~\ref{example alg ext}, we know that if $F$ has rank
    two, then $\langle x\rangle \leq F$ is algebraic if and only if
    $x$ is not a power of a primitive element of $F$.  Hence,
    situations like $H=\langle[a,b]\rangle < \langle a,b \rangle$ are
    examples of algebraic extensions where the base group $H$ is
    e-algebraically closed.  This is a behavior significantly
    different from what happens in field theory.
\end{example}

\begin{cor}\label{decide ealg closed}
Let $H\le F(A)$ be an extension of free groups.  If $H$ is finitely
generated, it is decidable whether $H$ is e-algebraically closed.
\end{cor}

\proof Let $x\not\in H$, viewed as a reduced word on the alphabet $A$,
let $p$ be the longest prefix of $x$ labeling a path starting at the
designated vertex 1 in $\Gamma_{A}(H)$, and let $s$ be the longest
suffix of $x$ labeling a path to 1 in $\Gamma_{A}(H)$.  We denote by
$1\cdot p$ and $1\cdot s\inv$ the end vertices of these two paths.

First assume that the sum of the length of $p$ and $s$ is less than
the length of $x$, that is, if $x = pys$ for some non-empty word $y$.
Then $\Gamma_{A}(\langle H,x\rangle)$ is obtained from $\Gamma_{A}(H)$
by gluing a path (made of new vertices and new edges) from $1\cdot p$
to $1\cdot s\inv$, labeled $y$.  In particular, $\rank(\langle
H,x\rangle) = \rank(H)+1$.

We now assume that the sum of the lengths of $p$ and $s$ is greater
than or equal to the length of $x$, and we let $t$ be the longest
suffix of $p$ which is also a prefix of $s$.  That is, $p = p't$, $s =
ts'$ and $x = p'ts'$.  Let $1\cdot p'$ be the end vertex of the path
starting at 1 and labeled $p'$ in $\Gamma_{A}(H)$.  If $1\cdot p' =
1\cdot s\inv$, then $x = p's$ labels in fact a loop at 1, that is,
$x\in H$, a contradiction.  So the labeled graph $\Gamma_{A}(\langle
H,x\rangle)$ is the quotient of $\Gamma_{A}(H)$ by the congruence
generated by the pair $(1\cdot p', 1\cdot s\inv)$ (see the end of
Section~\ref{sub of sub}).

Thus, in view of Corollary~\ref{charact ealg closed}, $H$ is
e-algebraically closed if and only if the following holds: for each
pair of distinct vertices $(v,w)$ in $\Gamma_{A}(H)$, the subgroup
represented by the quotient of $\Gamma_{A}(H)$ by the congruence
generated by $(v,w)$ has rank at most $\rank(H)$.  This is decidable,
and concludes the proof.  \eop

\section{Abstract properties of subgroups}\label{sec abstract}

Let $F$ be a free group.  An abstract \emph{property} of subgroups of
$F$ is a set $\p$ of subgroups of $F$ containing at least the total
group $F$ itself.  For simplicity, if $H\in \p$, we
will say that the subgroup $H$ \emph{satisfies property $\p$}.

We say that the property $\p$ is \emph{(finite) intersection closed}
if the intersection of any (finite) family of subgroups of $F$
satisfying $\p$ also satisfies $\p$, and that it is \emph{free factor
closed} if every free factor of a subgroup of $F$ satisfying $\p$ also
satisfies $\p$.  Finally, we say that the property $\p$ is
\emph{decidable} if there exists an algorithm to decide whether a
given finitely generated subgroup $H\leq F(A)$ satisfies $\p$.

\subsection{$\calP$-closure of a subgroup}

Let $F$ be a free group, $\p$ be an abstract property of subgroups of
$F$, and let $H\leq F$.  If there exists a unique minimal subgroup of
$F$ satisfying $\p$ and containing $H$, it is called the
\emph{$\p$-closure} of $H$, denoted by $\cl_{\p}(H)$; in this
situation, we say that $H$ admits a \textit{well defined}
$\p$-closure.

\begin{prop}\label{existenceofclosure}
Let $F$ be a free group and let $\p$ be an abstract property of
subgroups of $F$.
\begin{itemize}
\item[(i)] If $\p$ is intersection closed, then every subgroup $H\leq
F$ admits a well defined $\p$-closure.
\item[(ii)] If $\p$ is finite intersection closed and free factor
closed then every finitely generated subgroup $H\lefg F$ admits a well
defined $\p$-closure.
\item[(iii)] If $\p$-closures are well defined and $\p$ is free factor
closed, then for every subgroup $H\leq F$, we have $H\lalg
\cl_{\p}(H)$.  In particular, if $H$ is finitely generated, then so is
$\cl_{\p}(H)$.
\end{itemize}
\end{prop}

\demo Statement (i) is immediate: it suffices to consider the
intersection of all the extensions of $H$ satisfying $\p$ (there is at
least one, namely $F$ itself).

If $\p$ is only finite intersection closed, but is also free factor
closed, we use Theorem~\ref{existcl}: since every extension of a
finitely generated subgroup $H$ is a free multiple of an algebraic
extension of $H$, then every extension of $H$ in $\p$ contains an
algebraic extension of $H$ in $\p$.  It follows that the intersection
of all extensions of $H$ in $\p$ is equal to the intersection of the
algebraic extensions of $H$ in $\p$.  But the latter intersection is
finite, and hence it satisfies $\p$ as well, which concludes the proof
of (ii).

Finally, if $\p$ is free factor closed, then $H$ is not contained in
any proper free factor of its $\p$-closure, that is, $H\lealg
\cl_{\p}(H)$.  \qed

It would be interesting to produce an example of an abstract property
$\p$ that is closed under free factors and finite intersections, not
closed under intersections, and non-trivial for finitely generated
subgroups (note that the property \textit{to be finitely generated}
satisfies the required closure and non-closure properties, but it is
trivial for finitely generated subgroups).

\begin{remark}
It is well known that the property of being normal in $F$ is closed
under intersections and not under free factors, and that given a
subgroup $H\le F$, the normal closure of $H$ is well-defined, and is
not in general finitely generated, even if $H$ is.
\end{remark}

\begin{prop}\label{comput}
Let $\p$ be an abstract property for subgroups of $F(A)$ for which
$\p$-closures are well defined.  If $\p$-closures of finitely
generated subgroups of $F(A)$ are computable, then $\p$ is decidable.
The converse holds if, additionally, $\p$ is free factor closed.
\end{prop}

\demo Let us assume that $\p$-closures are computable.  Then, in order
to decide whether a given $H\lefg F(A)$ satisfies $\p$, it suffices to
compute $\cl_{\p}(H)$, and to verify whether $H = \cl_{\p}(H)$.

Conversely, suppose that $\p$ is free factor closed and decidable.
Then, given $H\lefg F(A)$, one can compute the set $\AE(H)$, check which
algebraic extensions of $H$ satisfy $\p$ and identify the minimal
one(s).  By Proposition~\ref{existenceofclosure}, only one of them is
minimal, and that one must be $\cl_{\p}(H)$.  \qed

\begin{remark}
Proposition~\ref{existenceofclosure} states that every property of
subgroups that is closed under (finite) intersections and under free
factors yields a well-defined closure operator for (finitely
generated) subgroups of $F$, that can be obtained by looking
exclusively at algebraic extensions.

A form of converse holds too: if $K\lefg F$, let $\p_{K}$ be the
following property.  A subgroup $L$ satisfies $\p_{K}$ if and only if
$L$ is a free factor of an extension of $K$.  Clearly, $F$ satisfies
this property, and one can verify that $\p_{K}$ is intersection and
free factor closed.  Moreover, one can use Proposition~\ref{existcl}
to verify that the $\p_{K}$-closure of a subgroup $H\le K$ is exactly
the $K$-algebraic closure of $H$.  In particular, for every algebraic
extension $H\lealg K$, $K$ is the $\calP$-closure of $H$ for some
well-chosen property $\calP$.
\end{remark}

\subsection{Some algebraic properties}\label{sec algebraic properties}

Let us recall the definition of certain properties of subgroups, that
have been discussed in the literature.  Let $H\le F$ be an extension
of free groups.  We say that $H$ is
\begin{itemize}
\item \emph{malnormal} if $H^g \cap H =1$ for all $g\in F\setminus H$;
\item \emph{pure} if $x^n\in H$, $n\ne 0$ implies $x\in H$ (this
property is also called being closed under radical, or being
isolated);
\item \emph{$p$-pure} (for a prime $p$) if $x^n \in H$, $(n,p) = 1$
implies $x\in H$;
\end{itemize}

The following results on malnormal and pure closure were first shown
in \cite[Section 13]{KM}.  The proof given here, while not
fundamentally different, is simpler and more general.
Corollary~\ref{rk pure goes down} below gives further properties of
these closures.

\begin{prop}\label{retract closure}
Let $F(A)$ be a free group.  The properties (of subgroups) defined by
malnormal, pure, $p$-pure ($p$ a prime), retract and e-algebraically
closed subgroups are intersection and free factor closed, and
decidable for finitely generated subgroups.

For each of these properties $\p$, each subgroup $H\le F(A)$ admits a
well-defined $\p$-closure $\cl_{\p}(H)$, which is an algebraic
extension of $H$.  Finally, if $H\lefg F(A)$, the $\p$-closure of $H$ has
finite rank and is computable.
\end{prop}

\demo The closure under intersections and free factors of malnormality
is immediate from the definition.  The decidability of malnormality
was established in \cite{BMR}, with a simple algorithm given in
\cite[Corollary 9.11]{KM}.

The closure under intersections of the properties of purity and
$p$-purity is immediate.  Now, assume that $K$ is pure, $H\leff K$,
and let $x$ be such that $x^n \in H$ with $n\ne 0$.  Since $K$ is
pure, we have $x\in K$, and we simply need to show that a free factor
of a free group $F$ is pure, which was established in
Example~\ref{example algebraic elements} above.  Thus purity is free
factor closed.  The proof of the same property for $p$-purity is
identical.  The decidability of purity and $p$-purity was proved in
\cite{BMMW94,BMMW}.

It is shown in~\cite[Lemma 18]{Be} that an arbitrary intersection of
retracts of $F$ is again a retract of $F$.  Moreover, it follows from
the definition of retracts that a retract of a retract is a retract,
and that a free factor is a retract.  Thus the property of being a
retract of $F$ is free factor closed.  The decidability of this
property was established by Turner, but as no proof seems to have
been published, we give his in Proposition~\ref{decide retracts}
below.

Suppose that $H\leff K\le F$, $K$ is e-algebraically closed and
$x\not\in H$.  If $x\not\in K$, then $\langle K,x\rangle = K*\langle
x\rangle$, so $\langle H,x\rangle =H*\langle x\rangle$.  If $x\in K
\setminus H$, then we have that $H$ is a free factor of $\langle H,x
\rangle \leq K$ and so, by Proposition~\ref{elext}, we also conclude
that $\langle H,x\rangle = H*\langle x\rangle$.  Thus the property of
being e-algebraically closed is closed under free factors.  Next, let
$(H_{i})_{i\in I}$ be a family of e-algebraically closed subgroups,
let $H = \bigcap_{i} H_{i}$, and let $x\not\in H$.  There exists $i\in
I$ such that $x\not\in H_{i}$, so $\langle H_{i},x\rangle = H_{i} *
\langle x\rangle$.  Using Lemma~\ref{technical e-algebraic}, we
conclude that $\langle H,x\rangle = H * \langle x\rangle$.  Thus the
property of being e-algebraically closed is also closed under
intersections.  Finally, this property is decidable by
Corollary~\ref{decide ealg closed}.

The last part of the statement follows from Proposition
\ref{existenceofclosure}.  \qed

As announced in the proof of Proposition~\ref{retract closure}, we
prove the decidability of retracts, that was established by
Turner~\cite{Ted private}.

\begin{prop}\label{decide retracts}
Let $H\le F(A)$ be an extension of finitely generated free groups.  It is
decidable whether $H$ is a retract of $F(A)$.
\end{prop}

\demo (Turner) Suppose that $A = \{a_{1},\ldots,a_{n}\}$ and let
$u_{1},\ldots, u_{r}$ be a basis of $H$.  Then $H$ is a retract of $F(A)$
if and only if there exist $x_{1},\ldots,x_{n}\in H$ such that the
endomorphism $\phi$ of $F(A)$ defined by $\phi(a_{i}) = x_{i}$ maps
$H$ identically to itself.  That is, if $u_{i}(x_{1},\ldots,x_{n}) =
u_{i}$ for $i = 1,\ldots,r$.  This can be expressed in terms of
systems of equations.

Let $e_{i}$ be the word on alphabet $\{X_{1},\ldots,X_{n}\}$ obtained
from the word $u_{i}$ (on alphabet $A$) by substituting $X_{j}$ for
$a_{j}$ for each $j$.  Then $H$ is a retract of $F$ if and only if the
system of equations $e_{i}(X_1,\ldots ,X_n) =u_i$, $i=1,\ldots ,r$
(where $u_i$ are viewed as constants in $H$) admits a solution in $H$.
This is decidable by Makanin's algorithm \cite{Makanin} (note that the
form of the system (i.e. the words $e_{i}$) depends on the way $H$ is
embedded in $F$, but once this form is established, the system itself
is entirely set within $H$, so Makanin's algorithm works, applied to
this system over $H$).  \qed

Let $H\le F$ be an extension of free groups.  Recall that $H$ is
\emph{compressed} if $\rank(H)\leq \rank(K)$ for every $K\leq F$
containing $H$ (see Section~\ref{sec alg closure}), and say that $H$
is \emph{inert} if $\rank(H\cap K)\leq \rank(K)$ for every $K\leq F$.
Both these properties were introduced by Dicks and Ventura \cite{DV}
in the context of the study of subgroups of free groups that are fixed
by sets of endomorphisms or automorphisms (see also~\cite{V2}).

It is clear that an inert or compressed subgroup is finitely
generated, with rank at most $\rank(F)$.  It is also clear that inert
subgroups (and retracts) of $F$ are compressed.  On the other hand, we
do not know whether all compressed subgroups are inert, nor whether
retracts are inert (both these facts are conjectured in~\cite{V2} and
related to other conjectures about fixed subgroups in free groups).

\begin{prop}\label{inert is closed}
Let $F$ be a free group.  The properties of inertness and
compressedness are closed under free factors.  In addition, inertness
is closed under intersections.

Each subgroup $H\le F$ admits an inert closure, which is an algebraic
extension of $H$.
\end{prop}

\demo The closure of inertness under intersections is shown
in~\cite[Corollary~I.4.13]{DV}.  Free factors of $F$ are trivially
inert.  Moreover, if $H\leq K\leq F$, $H$ is inert in $K$ and $K$ is
inert in $F$, then $H$ is inert in $F$.  So inertness is also closed
under free factors.

Now suppose that $H = L * M \leq F$ is compressed, and let $L\leq
K\leq F$.  Since $H\leq \langle K, M \rangle$, we have $$
\rank(L) + \rank(M) = \rank(H) \leq \rank(\langle K, M \rangle)\leq
\rank(K) + \rank(M).  $$
It follows that $\rank(L)\leq \rank(K)$, and hence $L$ is compressed.
Thus, compressedness is closed under free factors.  The last statement
is a direct application of Proposition~\ref{existenceofclosure}.  \qed

Note that, even though a finitely generated subgroup $H$ admits an
inert closure, which is one of its (finitely many) algebraic
extensions of $H$, we do not know how to compute this closure, nor how
to decide whether a subgroup is inert.

It is not known either whether compressedness is closed under
intersections, or even finite intersections, so we don't know whether
each subgroup admits a \textit{compressed closure}.  However it is
decidable whether a finitely generated subgroup of $F$ is compressed
\cite{V}.  Indeed if $H\lefg F$, then $H$ is compressed if and only if
$\rank(H)\leq \rank(K)$ for every algebraic extension $H\lalg K\leq
F$, which reduces the verification to a finite number of rank
comparisons.

\subsection{On certain topological closures}\label{sec topologies}

Let $\calT$ be a topology on a free group $F$.  The abstract property
of subgroups consisting of the subgroups that are closed in $\calT$ is
trivially closed under intersections.  This property becomes more
interesting when the topology is related to the algebraic structure of
$F$.  This is the case of the pro-\V\ topologies that we analyze now.

A \textit{pseudovariety of groups} \V\ is a class of finite groups
that is closed under taking subgroups, quotients and finite direct
products.  \V\ is called \textit{non-trivial} if it contains some
non-trivial finite group.  Additionally, if for every short exact
sequence of finite groups, $1\to G_1\to G_2\to G_3\to 1$, with $G_1$
and $G_3$ in \V, one always has $G_2\in \textbf{V}$, we say that \V\
is \textit{extension-closed}.

For every non-trivial pseudovariety of groups \V, the \textit{pro-\V}
topology on a free group $F$ is the initial topology of the collection
of morphisms from $F$ into groups in \V, or equivalently, the topology
for which the normal subgroups $N$ such that $F/N\in\V$ form a basis
of neighborhoods of the unit.  We refer the readers to
\cite{MSW,Weil00} for a survey of results concerning these topologies
with regard to finitely generated subgroups of free groups.  In
particular, Ribes and Zalesski\u\i\ showed that if \V\ is
extension-closed then every free factor of a closed subgroup is closed
\cite{RZ94}.  The following observation then follows from
Proposition~\ref{existenceofclosure}.

\begin{obs}\label{fact topological closure}
Let \V\ be a non-trivial extension-closed pseudovariety of groups.
Then the pro-\V\ closure of a finitely generated subgroup $H$ is
finitely generated, and an algebraic extension of $H$.
\end{obs}

In the case of the pro-$p$ topology ($p$ is a prime and the
pseudovariety \V\ is that of finite $p$-groups, which is closed under
extensions), Ribes and Zalesski\u\i\ \cite{RZ94} showed that one can
compute the closure of a given finitely generated subgroup of $F(A)$.  A
polynomial time algorithm was later given by Margolis, Sapir and Weil
\cite{MSW}, based on the finiteness of the number of principal
overgroups of $H$, that is, essentially on the spirit of
Fact~\ref{fact topological closure}.  Moreover, they showed that one
can simultaneously compute the pro-$p$ closures of $H$, for all primes
$p$, using the fact that they are all algebraic extensions, and hence
that they take only finitely many values.  This was also used to show
the computability of the pro-nilpotent closure of a finitely generated
subgroup: even though the pseudovariety of finite nilpotent groups is
not closed under extensions, it still holds that the pro-nilpotent
closure of a finitely generated subgroup is finitely generated and
computable.

At this point, several remarks are in order.  First, Ribes and
Zalesski\u\i\ \cite{RZ94} proved that if \V\ is extension-closed and
if $\bar{H}$ is the pro-\V\ closure of $H$, then
$\rank(\bar{H})\le\rank(H)$.  The proof of this fact can be reduced to
dimension considerations in appropriate vector spaces.  This proof
does not seem related with the idea of e-algebraic extensions, which
also lowers the rank (Corollary~\ref{rk goes down if elem*}).

Next, not every algebraic extension arises as a pro-\V\ closure for
some \V. This is clear if $H\lealg K$ and $\rank(K) > \rank(H)$ by the
result of Ribes and Zalesski\u\i\ cited above, but rank is not the
only obstacle.  Consider indeed $H = \langle a, bab\inv\rangle \le
F(a,b)$.  Then $H \lealg F$ (Example~\ref{example alg ext}) and
$\AE(H) = \calO_{A}(H) = \{H,F\}$.  We now verify that $H$ is
\V-closed for each non-trivial extension-closed pseudovariety \V, so
$F$ is never the \V-closure of $H$.  Since \V\ is non-trivial, the
cyclic $p$-element group $C_{p} = \langle c \mid c^p\rangle$ sits in
\V\ for some prime $p$.  Let $\phi_{p}\colon F \rightarrow C_{p}$ be
the morphism defined by $\phi_{p}(a) = 1$ and $\phi_{p}(b) = c$, and
let $N_{p} = \ker\phi_{p}$.  Then $H\le N_{p}$ and $N_{p}$ is
\V-closed, so $H$ is not topologically dense in $F$.  Since the
\V-closure of $H$ is in $\AE(H)$, it follows that $H$ is closed in the
pro-\V\ topology.

Solvable groups form an extension-closed pseudovariety, so the above
results apply to it: in particular, given an extension $H\lefg F(A)$, we
can compute a finite list of candidates for being the pro-solvable
closure of $H$, namely $\AE(H)$ (or even this list, restricted to the
extensions of rank at most $\rank(H)$).  However, it is a wide open
problem to compute this closure.

Finally, let us consider the (uncountable) collection of
extension-closed pseudovarieties of finite groups \V\ as above.  For
each finitely generated subgroup $H\le F$, the pro-\V\ closures of $H$
are among the (finitely many) algebraic extensions of $H$, so each
finitely generated subgroup $H$ naturally induces a finite index
equivalence relation on the collection of the \V's.  It would be
interesting to investigate the properties of these equivalence
relations.  In particular, the intersection of these equivalence
relations, as $H$ runs over all the (countably many) finitely
generated subgroups of $F(a,b)$, has countably many classes, so there
are pseudovarieties \V\ that are indistinguishable in this way.

\subsection{Equations over a subgroup}\label{sec equations}

In this section we use equations over free groups to define abstract
properties of subgroups.  Let $H\leq F$ be an extension of free
groups.  A \emph{(one variable) $H$-equation} (or \textit{equation
over $H$}) is an element $e = e(X)$ of the free group $H*\langle X
\rangle$, where $X$ is a new free letter, called the \emph{variable}.
An element $x\in F$ is a \emph{solution of $e(X)$} if $e(x) = 1$ in
$F$ (technically: if the morphism $H * \langle X \rangle \rightarrow
F$ mapping $H$ identically to itself and $X$ to $x$, maps $e$ to $1$).

\begin{example}
If $H = \langle a^2 \rangle$, the $H$-equation $e(X) = Xa^2
X^{-1}a^{-2}$ admits $a$ as a solution.  So does the $H$-equation
$X^2a^{-2}$.

If $e$ does not involve $X$, that is, $e\in H$, then $e$ has no
solution unless it is the \emph{trivial} equation $e = 1$, in which
case every element of $F$ is a solution.
\end{example}

We immediately observe the following.

\begin{lem}\label{sol-alg}
Let $H\le F$ be an extension of free groups and let $x\in F$.  The
element $x$ is a solution of some non-trivial $H$-equation if and only
if the elementary extension $H\le \langle H,x\rangle$ is algebraic.
\end{lem}

\demo Let $X$ be a new free generator and let $\phi \colon H * \langle
X \rangle \to F$ be the morphism that maps $H$ identically to itself
and $X$ to $x$.  By definition, $x$ is a solution of some non-trivial
equation over $H$ if and only if $\phi$ is not injective, and we
conclude by Proposition~\ref{elext} and Corollary~\ref{+1} that this
is equivalent to $H \lealg \langle H, x\rangle$.  \qed

In order to make this natural definition of equations independent on
the choice of the subgroup $H$, we consider a countable set
$X,Y_1,Y_2,\ldots$ of variables and we call \emph{equation} any
element $e$ of the free group on these variables.  If $H\le F$ is an
extension of free groups, a \emph{particularization of $e$ over $H$}
is the $H$-equation $e(X,h_1,h_2,\ldots)$ obtained by substituting
elements $h_1,h_2,\ldots \in H$ for the variables $Y_1,Y_2,\ldots$
(and having $X$ as variable).

A \emph{solution of the equation $e$ over $H$} is a solution of some
non-trivial particularization of $e$ over $H$, that is, an element
$x\in F$ such that, for some $h_1,h_2\ldots \in H$,
$e(X,h_1,h_2,\ldots) \neq 1$ but $e(x,h_1,h_2,\ldots) = 1$.  (Note
that even when $X$ occurs in $e$, some particularizations of $e$ over
$H$ can be trivial).

Let $\mathcal{E}$ be an arbitrary set of equations.  We say that a
subgroup $H\le F$ is \emph{$\calE$-closed} if $H$ contains every
solution over $H$ of every equation in $\mathcal{E}$.  Note that, when
looking for solutions, the set $\mathcal{E}$ is not considered as a
system of equations, but as a set of mutually unrelated equations.  In
particular, a larger set $\calE$ yields a larger set of solutions.

\begin{prop}\label{eqprop}
Let $F$ be a free group and let $\mathcal{E}$ be a set of equations.
Then the property of being $\mathcal{E}$-closed is closed under
intersections and under free factors.
\end{prop}

\demo The closure under intersections follows directly from the
definition.  Now assume that $K\le F$ is $\calE$-closed and let
$H\leff K$.  Let $x$ be a solution of an equation of $\calE$ over $H$.
Then $x$ is also a solution over $K$, and hence $x\in K$.  Now, by
Lemma~\ref{sol-alg}, $H \lealg \langle H, x\rangle \le K$.  This
contradicts $H\leff K$ unless $H = \langle H,x\rangle$, and hence
$x\in H$.  \qed

\begin{cor}\label{main on equations}
Let $H\le F$ and let $\calE$ be a set of equations.  There exists a
least $\calE$-closed extension of $H$, denoted by $\cl_{\calE}(H)$ and
called the $\calE$-closure of $H$.  Moreover, $H\lealg
\cl_{\calE}(H)$.

If in addition $H$ is finitely generated, then
$H\leealg\cl_{\calE}(H)$, $\rank(\cl_{\calE}(H)) \le \rank(H)$ and
there exists a finite subset $\calE_{0}$ of $\calE$ such that
$\cl_{\calE_{0}}(H) = \cl_{\calE}(H)$.
\end{cor}

\demo Propositions \ref{existenceofclosure} and \ref{eqprop} directly
prove the first part of the statement.

We now suppose that $H\lefg F$ and we let $H_{0} = H$ and suppose that
we have constructed distinct extensions $H_{0} \leealg H_{1} \leealg
\cdots \leealg H_{n}$ ($n\ge 0$), elements $x_{1},\ldots, x_{n}\in F$,
and equations $e_{1},\ldots,e_{n} \in \calE$ such that $H_{i} =
\langle H_{i-1}, x_{i}\rangle$ and $x_{i}$ is a solution of $e_{i}$
over $H_{i-1}$.  If $H_{n}$ is not $\calE$-closed, then there exists
an equation $e_{n+1}\in\calE$, and an element $x_{n+1}\not\in H_{n}$
such that $x_{n+1}$ is a solution of a non-trivial particularization
of $e_{n+1}$ over $H_{n}$.  Then $H_{n+1} = \langle H_{n},
x_{n+1}\rangle$ is a proper elementary algebraic extension of $H_{n}$
by Lemma~\ref{sol-alg}.  Since $H$ has only a finite number of
algebraic extensions, this construction must stop, that is, for some
$n$, $H_{n}$ is $\calE$-closed.  It follows easily that $H_{n}$ is the
$\calE$-closure of $H$, whose existence was already established.  In
particular $H\leealg\cl_{\calE}(H)$, and $\rank(\cl_{\calE}(H)) \le
\rank(H)$ by Corollary~\ref{rk goes down if elem*}.

Finally, let $\calE_{0} = \{e_{1},\ldots,e_{n}\}$.  Any $\calE$-closed
subgroup is also $\calE_{0}$-closed, and the $\calE_{0}$-closure of
$H$ must contain $H_{1},\ldots, H_{n}$.  Thus $\cl_{\calE}(H) =
\cl_{\calE_{0}}(H)$.  \qed

We conclude with the observation that some of the properties discussed
in Section~\ref{sec algebraic properties} can be expressed in terms of
equations.  Let $p$ be a prime number and let
\begin{eqnarray*}
\mathcal{E}_{mal} &=& \{X^{-1}Y_1 XY_2\},\cr
\mathcal{E}_p &=& \{X^n Y_1 \mid (n,p) = 1\},\cr
\mathcal{E}_{\mathbb{Z}} &=& \{X^n Y_1 \mid n\neq 0\}
=\bigcup_{p}\calE_{p},\cr
\mathcal{E}_{com} &=& \{X^{-1}Y_1^{-1} XY_1\}.
\end{eqnarray*}

\begin{prop}
Let $H\leq F$ be an extension of free groups. The subgroup $H$ is
\begin{itemize}
\item[(i)] malnormal if and only if it is $\calE_{mal}$-closed;
\item[(ii)] $p$-pure if and only if it is $\calE_{p}$-closed;
\item[(iii)] pure if and only if it is $\calE_{\mathbb{Z}}$-closed,
and if and only if it is $\calE_{com}$-closed.
\end{itemize}
\end{prop}

\demo $H$ is $\calE_{mal}$-closed if and only if, for all $h_1,h_2 \in
H$, not simultaneously trivial, every solution of the equation
$X^{-1}h_1 X h_2 = 1$ belongs to $H$.  That is, if and only if
$x^{-1}Hx \cap H\neq 1$ implies $x\in H$.  This is precisely the
malnormality property for $H$.  This proves~(i).

$H$ is $\calE_{p}$-closed if and only if $H$ contains the $n$-th roots
of every one of its elements, for all $n$ such that $(n,p) = 1$.
Again, this is exactly the definition of $p$-purity, showing~(ii).

Similarly, $H$ is $\calE_{\mathbb{Z}}$-closed if and only if $H$ is
pure.  Finally, we recall that two elements $x$ and $y$ in $F$ commute
if and only if they are powers of a common $z\in F$.  Thus the
subgroup generated by $H$ and all the roots of its elements is exactly
the $\calE_{com}$-closure of $H$.  \qed

Corollary \ref{main on equations} immediately implies the following.

\begin{cor}\label{rk pure goes down}
Let $H\lefg F$ and let $K$ be the malnormal (resp.  pure, $p$-pure)
closure of $H$.  Then $H\leealg K$ and $\rank(K) \le \rank(H)$.
\end{cor}

\section{Some open questions}\label{sec open}

To conclude this paper, we would like to draw the readers' attention
to a few of the questions it raises.
\begin{itemize}
    \item[(1)] We believe that the algebraic extensions of a finitely
    generated subgroup $H\lefg F$ are precisely the extensions which
    occur as principal overgroups of $H$ for every choice of an
    ambient basis.  That is, we conjecture that $\AE(H) = \bigcap_{A}
    \calO_{A}(H)$, where $A$ runs over all the bases of $F$.  As
    noticed in section~\ref{defs section}, this is the case when
    $H\lefi F$ or $H\leff F$, but nothing is known in general.
    
    \item[(2)] With reference to Corollary~\ref{rk goes down if
    elem*}, we would like to find an algebraic extension $H\lealg K$
    of finitely generated groups, where $\rank(K)\le \rank(H)$, yet
    the extension is not e-algebraic.  It would be appropriate to look
    for such an extension where $H$ is e-algebraically closed in $K$,
    that is, $\langle H,x\rangle = H*\langle x\rangle$ for each $x\in
    K\setminus H$ (Corollary~\ref{charact ealg closed}).
    
    \item[(3)] Even though a finitely generated subgroup $H$ admits an
    inert closure, which is one of the finitely many (computable)
    algebraic extensions of $H$, we do not know how to compute this
    closure.  Equivalently, it would be interesting to find an
    algorithm to decide whether a subgroup is inert (see
    section~\ref{sec algebraic properties}).
    
    \item[(4)] It is not known whether an intersection, even a finite
    intersection, of (strictly) compressed subgroups is again
    (strictly) compressed.  In other words, does a finitely generated
    subgroup admit a \textit{(strictly) compressed closure}?  If the
    answer was affirmative, then these closures would be computable,
    as indicated in section~\ref{sec algebraic properties}.
    
    \item[(5)] As pointed out in section~\ref{sec topologies}, we know
    that if \V\ is a non-trivial extension-closed pseudovariety of
    groups and $H\lefg F$, then $\bar H$, the pro-\V\ closure of $H$,
    is an algebraic extension of $H$ with rank at most $\rank(H)$.
    However the known proof of this fact does not rely on the notion
    of e-algebraic extensions.  We would like to find an example of
    such a subgroup $H$ and a pseudovariety \V\ such that the
    extension $H\le \bar H$ is not e-algebraic -- or alternately to
    give a new proof of Ribes and Zalesski\u\i's result (that in this
    situation, $\rank(\bar H)\le \rank(H)$), by showing that $H\leealg
    \bar H$.
    
    \item[(6)] As indicated at the end of section~\ref{sec
    topologies}, it would be interesting to find and investigate
    explicit examples of pseudovarieties $\text{\V }_1$ and $\text{\V
    }_2$, such that the pro-$\text{\V }_1$ and pro-$\text{\V }_2$
    closures of $H$ do coincide, for every $H\lefg F$.  As argued
    above, there are uncountably many such pairs being
    indistinguishable by means of closures of finitely generated
    subgroups.
    
    \item[(7)] Finally, Corollary~\ref{main on equations} shows that
    for every set of equations $\calE$ and every $H\lefg F$, there
    exists a finite subset $\calE_{0} \subseteq \calE$ such that
    $\cl_{\calE_{0}}(H)=\cl_{\calE}(H)$.  Is it true that such a
    finite set always exists satisfying the previous equality for all
    finitely generated subgroups of $F$ at the same time (showing a
    kind of noetherian behavior)?
\end{itemize}

\section*{Acknowledgements}

E. Ventura thanks the support received from DGI (Spanish government)
through grant BFM2003-06613, and from the Generalitat de Catalunya
through grant ACI-013.  P. Weil acknowledges support from the European
Science Foundation program \textsc{AutomathA}.  E. Ventura and P. Weil
wish to thank the Mathematics Department of the University of Nebraska
(Lincoln), where they were invited Professors when part of this work
was developped.  All three authors gratefully acknowledge the support
of the Centre de Recerca Matem\`atica (Barcelona) for their warm
hospitality during different periods of the academic year 2004-2005,
while part of this paper was written.


\end{document}